\documentclass[reqno]{amsart}

\usepackage[english]{babel}
\usepackage{graphicx}
\usepackage{textcomp}
\usepackage{amsmath}

\def\Div{\mathop{\rm div}\nolimits}

\usepackage{amssymb}
\usepackage{latexsym}
\newtheorem{proposition}{Proposition}
\newtheorem{definition}{Definition}
\newtheorem{remark}{Remark}
\newtheorem{Theorem}{Theorem}
\title[Probabilistic model associated with the pressureless gas dynamics]
{Probabilistic model associated with the pressureless gas dynamics}

\author[Albeverio,\,Korshunova,\, Rozanova]{Sergio Albeverio $^{1}$, Anastasia Korshunova $^{2}$,\, Olga Rozanova $^{2}$}

\address[$^{1}$]{Universit\"{a}t Bonn,
Institut f\"{u}r Angewandte Mathematik, Abteilung f\"{u}r
Stochastik, Endenicher Allee 60, D-53115 Bonn; HCM, SFB611 and IZKS,
Bonn; BiBoS, Bielefeld--Bonn; CERFIM, Locarno}
\address[$^{2}$]{Mathematics and Mechanics Faculty, Moscow State University, Moscow
119992, Russia}

\thanks {Supported by  DFG 436 RUS 113/823/0-1 (S.A. and O.R.) and the special program of the
Ministry of Education of the Russian Federation "The development of
scientific potential of the Higher School", project 2.1.1/1399 (A.K.
and O.R.).}

\email[$^{1}$]{albeverio@uni-bonn.de}
\email[$^{2}$]{rozanova@mech.math.msu.su}


\begin {document}

\maketitle

%\index{Albeverio S}
%\index{Korshunova A.A.}                              % write this for each author
%\index{Rozanova O.S.}                              % to generate the index
                            %

\begin{abstract}
Using a method of stochastic perturbation of a Langevin system
associated with the non-viscous Burgers equation we construct a
solution to the Riemann problem for the pressureless gas dynamics
describing sticky particles. As a bridging step we consider a medium
consisting of noninteracting particles. We analyze the difference in
the behavior of discontinuous solutions for these two models and the
relations between them. In our framework we obtain a unique entropy
solution to the Riemann problem in 1D case. Moreover, we describe
how starting from smooth data a $\delta$ - singularity arises in one
component of the solution.
\end{abstract}

\section*{Introduction}

We propose a method for solving the Riemann problem as well as for
describing the formation of singularities for the pressureless gas
dynamics system and a natural extension of it. The system of
pressureless gas dynamics is  very important since it is believed to
be the simplest model describing the formation of structures in the
universe (e.g.\cite{Shand}) and plays a significant role in the
theory of cooling gases and granular materials \cite{Brilliantov}.
It is a system consisting of two equations for the components of the
density $f$ and velocity $u$ expressing  the conservation of mass
and momentum \begin{equation}\label{sist_pred1}\partial_t
f+\Div_x(fu)=0,\end{equation}
\begin{equation}\label{sist_pred2}
\partial_t (fu)+\nabla_x(fu \otimes u)=0.\end{equation}
%(see below (\ref{sist_pred1}), (\ref{sist_pred2})).
 It first appears
to be very simple, however a closer analysis reveals that it has
some peculiar features due to its non strict hyperbolicity. The
system has attracted a significant interest in the last decades and
has been investigated quite intensively. In particular, it is well
known that the arising  in the velocity component of unbounded space
derivatives implies the generation of a $\delta$ - singularity in
the component of the density. Therefore for this system one needs to
define a generalized or measure-valued solution of a special kind.
This was done  in \cite{ERykovSinai}, \cite{brenier},
\cite{shelkv},\cite{DanShelkvDE},\cite{ShengZhang}, \cite{Joseph},
\cite{bouchut}, where the authors used different techniques
(vanishing viscosity, weak asymptotics, variational principle,
duality) to define the solution and prove its existence. Further,
the Riemann problem for the pressureless gas dynamics was studied
(e.g. \cite{ShengZhang}, \cite{shelkv},\cite{shelkvFPM}), including
a singular Riemann problem with a $\delta$ - singularity
concentrated at the jump at the initial moment. Nevertheless, even
in the 1D case there are certain open problems, not to mention those
present in the higher dimensional situation. In particular, there is
a problem concerning the uniqueness of solutions. Both in the case
of rarefaction and contraction it is possible to construct a whole
family of solutions to the Riemann problem satisfying the integral
identities and entropy conditions that are used to single out the
unique solution in the strictly hyperbolic case (see
\cite{danilov},\cite{huang} for details). Further, the process of
singularity formation was described up to now only in a very special
situation (\cite{DanilovMitrovic1},\cite{DanilovMitrovic2}).

Our method allows to find a unique solution to the Riemann problem
with arbitrary smooth left and right states. The most clear and
explicit results we get concern the case of constant left and right
states. We restrict ourselves to the 1D case, however the formulas
that we use are written in a similar way in any dimensions, and our
technique can be straightforwardly extended to the case of higher
dimensions. Of course the situation in higher dimension is much more
complicated, nevertheless, there are no fundamental obstacles to be
faced when applying our method.

Further, our method is constructive. It allows to describe the
behavior of the system starting from any smooth initial data. In
particular, it is possible to describe the singularity formation
including the time, position  and value of the amplitude of the
$\delta$ - function in the component of the density basing on the
initial data. Moreover, it is possible to describe the behavior of a
solution after the critical time, in particular, the position and
the amplitude of the $\delta$ - singularity.

%(\ref{sist_pred1}), (\ref{sist_pred2}) was considered.

Let us describe shortly the method in the 1D case. First of all it
is evident that till the solution to the system of pressureless gas
dynamics keeps its $C^1$ - smoothness, the velocity satisfies  the
non-viscous Burgers (or free transport) equation. We write it in the
Langevin form and introduce a  stochastic perturbation along the
trajectories of the particles. Further, we consider the position $x$
and velocity $u$ of the particles as random variables and find the
common probability density $P(t,x,u)$ in the space of positions and
velocities as a solution to the corresponding Fokker-Planck
equation. Then we introduce a pair of variables: the density of
particles $\rho(t,x)$ and the conditional expectation  of velocities
$\hat u(t,x)$  at fixed coordinates (and time).

We define a generalized solution in the sense of free particles
(FP), a pair \\$((f_{FP}(t,x), u_{FP}(t,x)),$ as a special double
limit of $(\rho,\hat u)$ as  the parameter $\sigma$ of stochastic
perturbation
 and the parameter $\varepsilon$ of the approximation
 of initial data  go to zero.

We prove that the pair $(\rho,\hat u)$  satisfies a gas dynamics
system with viscous and integral terms. The viscous term is the
usual one for the viscous approximation of the solution and it
vanishes as the parameter of stochastic perturbation $\sigma$ goes
to zero. The integral term vanishes if  $u_{FP},$ a part of the FP -
generalized solution,   is continuous and persists otherwise. Thus,
the FP - generalized solution to the pressureless gas dynamics
system solves in fact the extended gas dynamics system with an
integral term in the usual sense of integral identities. This
integral term can be considered as a spurious pressure, it is equal
to the dispersion of $u$ with respect to $P_x(t,x,u)$ ($P_x$
denoting the derivation with respect to $x$).

The FP - solution corresponds to the model describing the behavior
of a medium consisting in a micro-level of non-interacting
particles. In the component of density of the FP-solution the
$\delta$ - singularity can arise only from the domain where  of the
 initial velocity  has a special form (see Sec.6).
 The
$\delta$ - singularity does not arise in the FP- solution to the
Riemann problem with constant left and right states provided the
$\delta$ - singularity was not concentrated at the jump initially:
instead of the $\delta$ - singularity we have an overlapping domain
of a non-zero measure (a spurious pressure arises in this
overlapping domain). The FP-solution is interesting in itself,
moreover, it can be used to construct a solution to the sticky
particles models, where the particles are assumed to move together
when they meet (the name of  sticky particles model is used as a
synonym of the pressureless gas dynamics).  We propose a method
based on the conservation of the mass and momentum that allows to
reduce the FP-solution to the solution of the sticky particles model
and find the position and the amplitude of the $\delta$ -
singularity in the density component. The spurious pressure
degenerates in the process of above reduction. It fact, the problem
of collapsing the overlapping domain into a point only arises for
the compression waves, since for the rarefaction waves the initial
jump in the velocity decays into a smooth profile, such that the
FP-solution and the solution in the sense of integral identities
coincide.

The paper is organized as follows. In Sec.1 we consider the model of
motion for free particles perturbed along their trajectories and
introduce the integral characteristics of the medium consisting of
these particles, in particular, the mean velocity at a fixed
coordinate $\hat u.$ In Sec.2 we study the properties of $\hat u$
and prove that the limit of this value as the parameter of the
stochastic perturbation go to zero (provided it is smooth) takes
part of solution to the pressureless gas dynamics system. In Sec.3
we give a notion of generalized solution in the sense of free
particles (FP) to the Cauchy problem for the pressureless gas
dynamics system. In Secs.4 and 5 in 1D case we find the FP-solution
for the classical and singular Riemann problems, respectively. In
Sec.6 we describe the arising of singularity from smooth initial
data for the pressureless gas dynamics system. In Sec.7 we discuss
the difference between the FP-solution and the generalized solution
in the sense of integral identities. In Sec.8 we propose a method of
changing the FP-solution to the solution in the sense of integral
identities, and discuss the solution to the Riemann problem with
non-constant states and the evolution of the singularity arising
from smooth data. In Sec.9 we extend the method to a class of
systems that can be obtained from the scalar conservation law with a
convex flux. Sec.10 is a conclusion where we discuss the related
approaches and methods.

\section{Stochastic perturbation of the Burgers equation}

Let us consider the Cauchy problem for the non-viscous Burgers
equation:
\begin{equation}
\label{equ_Burg} \partial_t u+(u,\nabla)u=0,\,t>0,\qquad
u(x,0)=u_0(x)\in C^1(\mathbb{R}^n)\cap C_b(\mathbb{R}^n).
\end{equation}
%subject to initial data $u(x,0)=u_0(x)$,
Here $u(x,t)=(u_1,...,u_n)(x,t)\,$ is a vector-function
$\mathbb{R}^{n+1}\rightarrow\mathbb{R}^n.$

It is well known that solving  this equation in the smooth setting
is equivalent to solving the system of ODEs
\begin{equation}
\label{ODU} \dot{x}(t)=u(t,x(t)),\,t>0,\quad \dot{u}(t,x(t))=0
\end{equation}
for the characteristics $x=x(t)$.

We associate with (\ref{ODU}) the following  system of stochastic
differential equations:
$$dX_k(t)=U_k(t)dt+\sigma d(W_k)_t,$$
\begin{equation}
\label{SDU} dU_k(t)=0,\quad k=1,...,n,\,
\end{equation}
$$X(0)=x,\quad U(0)=u,$$
where $X(t)$ and $U(t)$ are considered as  random variables with
given initial distributions, $(X(t),U(t))$ runs in the phase space
$\mathbb{R}^n\times\mathbb{R}^n,$  $\sigma$ is a real strictly
positive constant and $(W)_t=(W_k)_t$, $k=1,...,n$ is an n -
dimensional Brownian motion, $\,t>0.$

System (3) describes the free motion of a particle that "does not
feel" the other particles. We assume that initially at time $t=0$
these non-interacting particles are distributed with a density
$f(x)$ and denote by $P(t,x,u)$  the probability density in position
and velocity space for the solutions of (\ref{SDU}),
$x\in\mathbb{R}^n,\, u\in\mathbb{R}^n$ at time $t.$

The stochastic system described by (\ref{SDU}) is associated with
the deterministic non-viscous Burgers equation (\ref{equ_Burg}) in
the sense that the deterministic characteristics (\ref{ODU}) are
replaced by a stochastic perturbation of the characteristics as
described by (\ref{SDU}). System (\ref{SDU}) with this
interpretation is what we understand as "stochastic perturbation of
the Burgers equation."

 Let us introduce the function
\begin{equation}
\label{u_sdu}
\hat{u}(t,x)=\dfrac{\int\limits_{\mathbb{R}^n}uP(t,x,u)du}{\int\limits_{\mathbb{R}^n}P(t,x,u)du},\quad
t\ge 0,\quad x\in \mathbb{R}^n.
\end{equation}
This value (\ref{u_sdu}) can be interpreted as the conditional
expectation of $U$ for fixed position $X$ \cite{Chorin}. If we
choose as initial distribution

\begin{equation}
\label{P0}
P_0(x,u)=\delta(u-u_0(x))f_0(x)=\prod\limits_{k=1}^n\delta(u_k-(u_0(x))_k)f_0(x),
\end{equation}
where $f_0$ is an arbitrary sufficiently regular nonnegative
function such that $\,\int\limits_{{\mathbb R}^n}f_0(x)dx=1$, then
$\hat{u}(0,x)=u_0(x)$. Certain properties of $\hat{u}(t,x)$ have
been established in \cite{Roz1}( see also \cite{Roz2} for another
type of stochastic perturbation).

The density $P=P(t,x,u)$ obeys the Fokker-Planck equation
\begin{equation}
\label{Fok-Plank} \dfrac{\partial P}{\partial
t}=\left[-\sum\limits_{k=1}^n u_k\dfrac{\partial}{\partial
x_k}+\sum\limits_{k=1}^n\dfrac12\sigma^2\dfrac{\partial^2}{\partial
x_k^2}\right] P,
\end{equation}
subject to the initial data (\ref{P0}).

We apply the Fourier transform to $P(t,x,u)$ in (\ref{Fok-Plank}),
(\ref{P0}) with respect to the  variables $x$ and $u$   and obtain
the Cauchy problem for the Fourier transform
$\tilde{P}=\tilde{P}(t,\lambda,\xi)$ of $P(t,x,u)$:
\begin{equation}
\label{preobr_Fok-Plank} \dfrac{\partial \tilde{P}}{\partial
t}=-\dfrac12\sigma^2|\lambda|^2\tilde{P}+\lambda\dfrac{\partial\tilde{P}}{\partial
\xi},
\end{equation}
\begin{equation}
\label{preobr_P0}\tilde{P}(0,\lambda,\xi)=\int\limits_{\mathbb{R}^n}e^{-i(\lambda,s)}e^{-i(\xi,u_0(s))}f_0(s)ds,
\qquad \lambda,\xi \in {\mathbb R}^n.
\end{equation}
Equation (\ref{preobr_Fok-Plank}) can easily be  integrated and we
obtain the solution given by the following formula:
\begin{equation}
\label{preobr_P}\tilde{P}(t,\lambda,\xi)=\tilde{P}(0,\lambda,\xi+\lambda
t)e^{-\frac12\sigma^2|\lambda|^2t}.
\end{equation}

The inverse Fourier transform (in the distributional sense) allows
to find the density $P(t,x,u),\,t>0$:
$$P(t,x,u)=\dfrac1{(2\pi)^{2n}}\int\limits_{\mathbb{R}^n}
\int\limits_{\mathbb{R}^n}e^{i(\lambda,x)}e^{i(\xi,u)}\tilde{P}(t,\lambda,\xi)\,d\lambda
d\xi=$$

$$=\dfrac1{(2\pi)^{2n}}\int\limits_{\mathbb{R}^n}\int\limits_{\mathbb{R}^n}e^{i(\lambda,x)}e^{i(\xi,u)}
\int\limits_{\mathbb{R}^n}e^{-i(\lambda,s)}e^{-i(\xi+\lambda
t,u_0(s))}f_0(s)ds e^{-\frac12\sigma^2|\lambda|^2t}d\lambda d\xi=$$

$$=\dfrac1{(2\pi)^{2n}}\int\limits_{\mathbb{R}^n}f_0(s)\int\limits_{\mathbb{R}^n}e^{i(\xi,u-u_0(s))}d\xi
\int\limits_{\mathbb{R}^n}e^{-\frac12\sigma^2t\left(\lambda-\frac{i|x-u_0(s)t-s|}{\sigma^2t}\right)^2-\frac{|u_0(s)t+s-x|^2}{2\sigma^2t}}d\lambda
ds=$$

\begin{equation}
\label{s_plotn}=\dfrac1{(\sqrt{2\pi
t}\sigma)^n}\int\limits_{\mathbb{R}^n}\,\delta(u-u_0(s))\,f_0(s)\,e^{-\frac{|u_0(s)t+s-x|^2}{2\sigma^2t}}ds,\quad
t\ge 0,\,x\in \mathbb{R}^n.
\end{equation}
Then we substitute $P(t,x,u)$ in (\ref{u_sdu}) and get the following
expression for $\hat{u}(t,x)$ (sometimes we insert a label $\sigma$
to stress the dependence on this parameter):
\begin{equation}
\label{sol_u_sdu}
\hat{u}(t,x)=\hat{u}_\sigma(t,x)=\dfrac{\int\limits_{\mathbb{R}^n}u_0(s)f_0(s)e^{-\frac{|u_0(s)t+s-x|^2}{2\sigma^2t}}ds}{\int\limits_{\mathbb{R}^n}f_0(s)e^{-\frac{|u_0(s)t+s-x|^2}{2\sigma^2t}}ds}.
\end{equation}

\begin{remark}
The integrals in (\ref{sol_u_sdu}) are defined also for a wider
class of $f_0$ than the probability density of the particle
positions in the space at the initial moment of time. If the
integral $\int\limits_{\mathbb{R}^n}f_0(x)dx$ diverges (for example,
for $f_0(x)=const$), we can consider the domain $[-L,L]^n$, where
$L>0$ and use another definition of $\hat{u}_\sigma(t,x)$:
\begin{equation}
\label{sol_u_sduL} \hat{u}_\sigma(t,x)=\lim\limits_{L\rightarrow
+\infty}\dfrac{\int\limits_{[-L,L]^n}u_0(s)f_0(s)e^{-\frac{|u_0(s)t+s-x|^2}{2\sigma^2t}}ds}
{\int\limits_{[-L,L]^n}f_0(s)e^{-\frac{|u_0(s)t+s-x|^2}{2\sigma^2t}}ds}
\end{equation}
(provided the limit exists). Evidently, this definition  coincides
with (\ref{sol_u_sdu}) for $f_0\in L_1(\mathbb{R}^n)$.
\end{remark}

\section{Properties of  velocity averaged at a fixed coordinate}

The following property of $\hat{u}(t,x)$ holds:

\begin{proposition}\label{prop1}
Let $u_0$ and $f_0>0 $ be  functions of class $C^1(\mathbb{R}^n)\cap
C_b(\mathbb{R}^n).$  If $ t_*(u_0)>0$ is a moment of time such that
the solution to the Cauchy problem (\ref{equ_Burg}) with the initial
condition $u_0$ keeps this smoothness for $0<t<t_*(u_0)\le+\infty,$
then $\hat{u}_\sigma(t,x)$ tends to a solution of problem
(\ref{equ_Burg}) as $\sigma\rightarrow 0$ for any fixed
$(t,x)\in\mathbb{R}^{n+1},\,0<t<t_*(u_0)$.
\end{proposition}

\proof
%\textsc{Proof}
Let us denote by $J(u_0(x))$ the Jacobian matrix of the map
$\,x\longmapsto u_0(x).$ As it was shown in \cite{protter}\,(Theorem
1), if $J(u_0(x))$ has at least one  eigenvalue which is negative
for a certain point $x\in{\mathbb R}^n,$  then the classical
solution to (\ref{equ_Burg}) fails to exist beyond a positive time
$t_*(u_0).$ Otherwise,  $t_*(u_0)=\infty.$ The matrix $C(t,x) = (I+t
J(u_0(x))),$ where $\,I\,$ is the identity matrix, fails to be
invertible for $t=t_*(u_0).$
%Further on we shall use the notation
%$(u_0)'_x$ instead of $|J(u_0(x)|.$

%The convergence for $\sigma\to 0$ is shown as follows.
The formula
(\ref{sol_u_sdu}) (or (\ref{sol_u_sduL})) implies, using the weak
convergence of measures and the fact that $f_0$ and $u_0$ are
continuous and bounded
$$\lim\limits_{\sigma\rightarrow
0}\hat{u}_\sigma(t,x)\,=\,\dfrac{\int\limits_{\mathbb{R}^n}u_0(s)f_0(s)\lim\limits_{\sigma\rightarrow
0}\frac1{(\sqrt{2\pi
t}\sigma)^n}e^{-\frac{|u_0(s)t+s-x|^2}{2\sigma^2t}}ds}{\int\limits_{\mathbb{R}^n}f_0(s)\lim\limits_{\sigma\rightarrow
0}\frac1{(\sqrt{2\pi
t}\sigma)^n}e^{-\frac{|u_0(s)t+s-x|^2}{2\sigma^2t}}ds}=$$
$$
\,\dfrac{\int\limits_{\mathbb{R}^n}u_0(s)f_0(s)\delta_{p(t,x,s)}ds}{\int\limits_{\mathbb{R}^n}f_0(s)
\delta_{p(t,x,s)}ds},
$$
with $p(t,x,s) = u_0(s)t+s-x\,,$ where $\delta_y$ is the Dirac
measure as $y\in {\mathbb R}^n.$ We can then on the basis on the
invertibility of $C(t,x)$ use locally the implicit function theorem
and find $s=s_{t,x}(p).$
%Moreover, $dp=\det C(t,s_{t,x})\,ds.$
Therefore,
$$\lim\limits_{\sigma\rightarrow
0}\hat{u}_\sigma(t,x)\,=\,
\dfrac{\int\limits_{\mathbb{R}^n}u_0(s_{t,x}(p))f_0(s_{t,x}(p))\,\det
(C(t,s_{t,x}(p)))^{-1}\,\delta_p\,(ds_{t,x})}
{\int\limits_{\mathbb{R}^n}f_0(s_{t,x}(p))\,\det(C(t,s_{t,x}(p)))^{-1}\,\delta_p\,(ds_{t,x})\,}=\,u_0(s_{t,x}(0)).
$$
Let us introduce the new notation  $s_0(t,x)\equiv s_{t,x}(0).$ Then
the following vectorial equation holds:
\begin{equation}
\label{usl}u_0(s_0(t,x))t+s_0(t,x)-x=0.
\end{equation}
%Thus,
%$$\lim\limits_{\sigma\rightarrow
%0}\hat{u}(t,x)=u_0(s_0(t,x)).$$
Let us show that $u(t,x)=u_0(s_0(t,x))$ satisfies the Burgers
equation, that is
\begin{equation}\label{Burgsubs}\sum\limits_{j=1}^n\,\partial_j
(u_{0,i})(s_{0,j})_t\,+\,\sum\limits_{j,k=1}^n  u_{0,j}
\partial_k(u_{0,i})(s_{0,k})_{x_j}=0, \quad i=1,...,n,
\end{equation}
and $u_0(s_0(0,x))=u_0(x)$. Here we denote by $u_{0,i}$ and
$s_{0,i}$ the $i$ - th components of vectors $u_0$ and $s_0,$
respectively.

We differentiate (\ref{usl}) with respect to $t$ and $x_j$ to get
the matrix equations: $$\sum\limits_{j=1}^n\,C_{ij}\,(s_{0,j})_t
+u_{0,i}=0,\, \quad i=1,...,n,$$ and
$$\sum\limits_{k=1}^n\,C_{ik}\,(s_{0,k})_{x_j}+\delta_{ij}=0,
\, \quad i,j=1,...,n,$$ where $\delta_{ij}$ is the Kronecker symbol.
The equations imply
\begin{equation}\label{subst}
(s_{0,j})_t\,=\,-\,\sum\limits_{i=1}^n\,(C^{-1})_{ij}\,u_{0,i},\qquad
(s_{0,k})_{x_j}\,=\,-\,(C^{-1})_{jk}.
\end{equation}
 It remains now only to
substitute (\ref{subst}) into (\ref{Burgsubs}).

Further,  (\ref{usl}) implies $u_0(s_0(0,x))=u_0(x)$.  $\,\square$

\medskip

It is important to note that $s_0(t,x)$ is unique for all $t$ for
which the solution to the Burgers equation $u(t,x)$ is smooth.
\medskip
\medskip
\begin{remark}Proposition \ref{prop1} can naturally be  extended to
the class of functions $f_0$ such that there exists a sequence
$f_0^\varepsilon\in C^1(\mathbb{R}^n)\cap C_b(\mathbb{R}^n)$
converging to $f_0$ as $\varepsilon\rightarrow 0$ almost everywhere.
In this case $\hat{u}_\sigma(t,x)$ tends to a solution of problem
(\ref{equ_Burg}) as $\sigma\rightarrow 0$ almost everywhere on
$(t,x)\in\mathbb{R}^{n+1},\,0<t<t_*(u_0)$.
%We will use this fact to study
%the Riemann problem solution in the case of rarefaction (see
%Sec.\ref{classRP}).\label{remf0_jump}

\end{remark}
\medskip
\medskip

Let us set $\rho(t,x)=\int\limits_{\mathbb{R}^n}P(t,x,u)du, \,t\ge
0,\, x\in {\mathbb R}^n$. From (\ref{s_plotn}) we have
\begin{equation}
\label{plotn} \rho(t,x)=\rho_\sigma(t,x)=\dfrac1{(\sqrt{2\pi
t}\sigma)^n}\int\limits_{\mathbb{R}^n}f_0(s)e^{-\frac{|u_0(s)t+s-x|^2}{2\sigma^2t}}ds,\,t\ge
0,\, x\in {\mathbb R}^n.
\end{equation}

\begin{proposition}
The scalar function $\rho(t,x)$ and the vector-function
$\hat{u}(t,x)=(\hat{u}_1,..,\hat{u}_n)\,$ defined in (\ref{u_sdu})
solve the following system:
\begin{equation}
\label{sist_obw1}\dfrac{\partial\rho}{\partial t}\,+\,\Div_x
(\rho\hat{u})\,=\,\dfrac12\sigma^2\sum\limits_{k=1}^{n}\dfrac{\partial^2\rho}{\partial
x_k^2},
\end{equation}
\begin{equation}
\label{sist_obw2}\dfrac{\partial(\rho\hat{u}_i)}{\partial t}\,+\,
\nabla(\rho\,\hat{u}_i\,\hat{u})\,=\,\dfrac12\sigma^2\sum\limits_{k=1}^{n}\dfrac{\partial^2(\rho\hat{u}_i)}{\partial
x_k^2}\,-\,\int\limits_{\mathbb{R}^n}(u_i\,-\,\hat{u}_i)\,\big((u\,-\,\hat{u}),\nabla_x
P(t,x,u)\big)\,du,
\end{equation}
$\quad i=1,..,n,\,t\ge 0.$
\end{proposition}

\proof

The equation (\ref{sist_obw1}) follows from the Fokker-Planck
equation (\ref{Fok-Plank}) directly.  To prove  (\ref{sist_obw2}) we
note that the definitions of $\hat{u}(t,x)$ and $\rho(t,x)$ imply
\begin{eqnarray}
%\begin{equation}
\label{p2_1}\dfrac{\partial(\rho\hat{u})}{\partial t}
=\dfrac{\partial}{\partial
t}\int\limits_{\mathbb{R}^n}uP(t,x,u)du=\int\limits_{\mathbb{R}^n}uP_t(t,x,u)\,du=\nonumber\\
=-\int\limits_{\mathbb{R}^n}u(u,\nabla_x
P(t,x,u))du+\dfrac12\sigma^2\sum\limits_{k=1}^{n}\dfrac{\partial^2\rho}{\partial
x_k^2},
%\end{equation}
\end{eqnarray}
where $P_t\equiv \frac{\partial }{\partial t} P.$

Further, for $i=1,..,n$ we have
\begin{eqnarray}\label{p2_2}
\dfrac{\partial (\rho\,\hat{u}_k\,\hat{u}_i)}{\partial x_k}=
\hat{u}_i\dfrac{\partial }{\partial
x_k}\left(\int\limits_{\mathbb{R}^n}\,u_k\,P\,du\right)\,+\,
\int\limits_{\mathbb{R}^n}\,u_k\,P\,du \,\dfrac{\partial }{\partial
x_k}\left(\frac{\int\limits_{\mathbb{R}^n}\,u_k\,P\,du}{\int\limits_{\mathbb{R}^n}\,P\,du}\right)=\nonumber\\
=\int\limits_{\mathbb{R}^n}\,\hat
u_i\,u_k\,P_{x_k}\,du\,+\,\int\limits_{\mathbb{R}^n}\,u_k\,P\,du\,
\frac{\int\limits_{\mathbb{R}^n}\,u_i\,P_{x_k}\,du
\,\int\limits_{\mathbb{R}^n}\,P\,du
\,-\,\int\limits_{\mathbb{R}^n}\,u_i\,P\,du\,\int\limits_{\mathbb{R}^n}\,P_{x_k}\,du
}{\left(\int\limits_{\mathbb{R}^n}\,P\,du\right) }=\\
=\int\limits_{\mathbb{R}^n}\,(u_k\,\hat u_i +u_i\,\hat u_k\,-\,\hat
u_k\,\hat u_i)\,P_{x_k}\,du,\quad i,k=1,...,n,\nonumber
\end{eqnarray}
with $P_{x_k} \,\equiv \,\frac{\partial }{\partial x_k} P.$

Equation (\ref{sist_obw2}) follows immediately from (\ref{p2_1}) and
(\ref{p2_2}). $\,\square$

\bigskip

In the one-dimensional case  (\ref{sist_obw1})-(\ref{sist_obw2}) has
the form
\begin{equation}\label{sist_obw1D1}\partial_t\rho+\partial_x
(\rho\hat{u})=\dfrac12\sigma^2\partial^2_{xx}\rho,\end{equation}
\begin{equation}\label{sist_obw1D2}
\partial_t(\rho\hat{u})+\partial_x
(\rho\hat{u}^2)=\dfrac12\sigma^2\partial^2_{xx}(\rho\hat{u})-\int\limits_{\mathbb{R}}(u-\hat{u})^2P_x(t,x,u)du.
\end{equation}

\bigskip

Let us set $f(t,x)=\lim\limits_{\sigma\rightarrow 0}\rho(t,x)$ and
$\bar u(t,x)=\lim\limits_{\sigma\rightarrow 0}\hat u(t,x).$

\begin{proposition}\label{prop3}
Assume that $(f(t,x),\bar u(t,x)),$ the limits  of $(\rho,\hat u)$
as $\sigma\to 0,$ are $C^1$ -- smooth bounded functions for
$(t,x)\in \Omega:=[0,t_*(u_0))\times {\mathbb R}^n,\, t_*(u_0)\le
\infty.$ Then they solve in $\Omega$ the pressureless gas dynamics
system (\ref{sist_pred1}), (\ref{sist_pred2}).
\end{proposition}

\proof As follows from Proposition \ref{prop1}, the function $\bar
u(t,x)$ is a $C^1$ - solution of the non-viscous Burgers equation.
Further,  (\ref{sist_obw1D1}) is a linear parabolic equation with
respect to $\rho,$ hence the limit as $\sigma\to 0$ reduces it to
the continuity equation (\ref{sist_pred1}). Equation
(\ref{sist_pred2}) is a corollary of the non-viscous Burgers
equation and (\ref{sist_pred1}). $\square$

%\begin{remark}

\begin{remark} Proposition \ref{prop3} implies that the integral
term on the right-hand side of (\ref{sist_obw2}) vanishes as
$\sigma\to 0$ in the case of  smooth limit functions $f$ and $\bar
u.$ Let us prove this fact alternatively. Indeed, as follows from
(\ref{s_plotn}), we have as $\sigma\to 0$
$$
\int\limits_{\mathbb{R}^n}(u_i\,-\,\hat{u}_i)\,\big((u\,-\,\hat{u}),\nabla_x
P(t,x,u)\big)\,du=$$
$$
\frac{1}{(\sqrt{2\pi t}\sigma)^n}\,\int\limits_{{\mathbb
R}^n}\,f_0(s)\,(u_{0i}(s)\,-\,\hat{u}_i)\,\big((u_0(s)\,-\,\hat{u})\,,\,\nabla_x
e^{-\frac{|u_0(s)t+s-x|^2}{2\sigma^2t}}\big)\,ds=
$$
$$
\frac{1}{(\sqrt{2\pi t}\sigma)^n t \sigma^2}\,\int\limits_{{\mathbb
R}^n}\,e^{-\frac{|u_0(s)t+s-x|^2}{2\sigma^2t}}\,f_0(s)\,(u_{0i}(s)\,-\,\hat{u}_i)\,\big((u_0(s)\,-\,\hat{u})\,,\,u_0(s)t+s-x
\big)\,ds=
$$
$$
\frac{1}{(\sqrt{2\pi t}\sigma)^n t
\sigma^2}\,\int\limits_{U_{s_0(t,x)}(\varepsilon)}\,e^{-\frac{|u_0(s)t+s-x|^2}{2\sigma^2t}}\,f_0(s)\,
(u_{0i}(s)\,-\,\hat{u}_i)\,\big((u_0(s)\,-\,\hat{u})\,,\,u_0(s)t+s-x
\big)\,ds\,+o(\sigma),
$$
where $U_{s_0(t,x)}(\varepsilon)$ is an $\varepsilon$ - neighborhood
of the point $s_0(t,x)$ (see the proof of Proposition \ref{prop1}),
$\,i=1,..,n$. Further,
$$|u_{0i}(s)\,-\,\hat{u}_i(t,x)|\,=\,|(u_{0i}(s)\,-\,\bar{u}_i(t,x))+ (\,\bar{u}_i(t,x)\,-\,\hat
u_{i}(t,x))\,|\,= $$
$$\,|(u_{0i}(s)\,-\,{\bar u}_{0i}(s_0(t,x)))+
(\,{\bar u}_{0i}(s_0(t,x))\,-\,\hat u_{i}(t,x))\,|.
$$
For every fixed $x$ and $t\in (0, t_*)$ and for every $\sigma>0$
there exists $\varepsilon(\sigma)>0$ such that if $s\in
U_{s_0(t,x)}(\varepsilon),$ then $|(u_{0i}(s)\,-\,{\bar
u}_{0i}(s_0(t,x)))<\sigma.$ Moreover, since $\hat
u_{\sigma_1}(t,x)\to \bar u(t,x)$ as $\sigma_1\to 0$ (we have
renamed the parameter), then for every $\sigma>0$ there exists
$\sigma_1(\sigma)>0$ such that for sufficiently small $\sigma_1 $ we
have $|\hat u_{i}(t,x)\,-\,{\bar u}_{0i}(s_0(t,x)|<\sigma.$

Thus,
$$
\big|\frac{1}{(\sqrt{2\pi t}\sigma)^n t
\sigma^2}\,\int\limits_{U_{s_0(t,x)}(\varepsilon)}\,e^{-\frac{|u_0(s)t+s-x|^2}{2\sigma^2t}}\,f_0(s)\,
(u_{0i}(s)\,-\,\hat{u}_i)\,\big((u_0(s)\,-\,\hat{u})\,,\,u_0(s)t+s-x
\big)\,ds\big|\le
$$
$$
{\rm Const}\,\cdot\,\frac{1}{(\sqrt{2\pi t}\sigma)^n}
\,\int\limits_{U_{s_0(t,x)}(\varepsilon)}\,e^{-\frac{|u_0(s)t+s-x|^2}{2\sigma^2t}}\,f_0(s)\,|u_0(s)t+s-x|\,
ds,
$$
where the constant does not depend of $\sigma.$
The latter integral tends to zero %in $\mathcal D'({\mathbb R}^n)$
as
$\sigma\to 0.$

In fact, to prove that the integral term vanishes as $\sigma\to 0,$
we have used only the continuity of $\bar u$ and the boundedness of
$f_0.$

However, as we will show in Sec.\ref{spurpress}, if  we put instead
of $\bar u$ a discontinuous function, this integral term does not
vanish.
\end{remark}

\medskip

%obtain in the limit $\sigma\rightarrow 0$ the system
%$\square$
\section{Generalized solution in the sense of free particles}

Being inspired by the fact  that the formula (\ref{sol_u_sdu}) makes
sense also for discontinuous initial data $(f_0(x),u_0(x)),$ we give
the following definition for any dimension.

\begin{definition}\label{definit1} We call the couple of functions $(f_{FP}(t,x),u_{FP}(t,x))$
{\it a generalized solution to the Cauchy problem
(\ref{sist_pred1}), (\ref{sist_pred2}) in the sense of free
particles (FP-generalized solution)}  subject to initial data $
(f_0(x),u_0(x)) \in {\mathbb L}^2_{loc}({\mathbb R}^n)\cap {\mathbb
L}_{\infty}({\mathbb R}^n)$, if for almost all $(t,x)\in
\mathbb{R}_+\times \mathbb{R}^n$
$$f_{FP}(t,x)=\lim\limits_{\varepsilon\rightarrow
0}(\lim\limits_{\sigma\rightarrow
0}\rho_\sigma^{\varepsilon}(t,x)),$$
$$u_{FP}(t,x)=\lim\limits_{\varepsilon\rightarrow
0}(\lim\limits_{\sigma\rightarrow
0}\hat{u}_\sigma^{\varepsilon}(t,x)),$$ where
$(\rho_\sigma^{\varepsilon}(t,x),\hat{u}_\sigma^{\varepsilon}(t,x))$
satisfy  the system (\ref{sist_obw1D1}), (\ref{sist_obw1D2}) with
initial data $$\bar f_0^{\varepsilon}\,=\, \eta_\varepsilon
*f_0,\quad \bar u_0^{\varepsilon}\,=\,\eta_\varepsilon
*u_0\,$$ where $\eta_\varepsilon(x)$ is the standard averaging
kernel.
\end{definition}

\begin{remark} The properties of the standard averaging
kernel \cite{Evans} imply that $f_0^{\varepsilon}$ and
$u_0^{\varepsilon}) $ belong to the class $C^\infty$ and
$$\lim\limits_{\varepsilon\rightarrow
0}\bar f_0^{\varepsilon}(x)=f_0(x),\quad
\lim\limits_{\varepsilon\rightarrow 0}\bar
u_0^{\varepsilon}(x)=u_0(x),$$ for almost all fixed
$x\in\mathbb{R}^n$.
\end{remark}

\begin{remark} As we will see below, if $u_{FP}$ is continuous,
the FP-generalized solution  is a solution to (\ref{sist_pred1}),
(\ref{sist_pred2}), for example, in the sense of integral
identities. However, if $u_{FP}$ is  discontinuous, the FP-solution
solves a different system, namely one  that differs from
(\ref{sist_pred1}), (\ref{sist_pred2}) by an integral term
corresponding to a spurious pressure. Nevertheless, using the
FP-solution we can solve (\ref{sist_pred1}), (\ref{sist_pred2})
itself.
\end{remark}

\begin{definition} We call the pair  $(f_0^{\varepsilon}, u_0^{\varepsilon})$
 a monotonic approximation of initial data $(f_0(x),u_0(x))\in
{\mathbb L}^2_{loc}({\mathbb R}^n)\cap {\mathbb L}_{\infty}({\mathbb
R}^n)$, if
\begin{itemize}
\item $f_0^{\varepsilon}$ and $u_0^{\varepsilon}$ are  from the class
$C_b(\mathbb{R}^n),$ moreover, they are from $C^1(\mathbb{R}^n)$
almost everywhere;
\item
$$\lim\limits_{\varepsilon\rightarrow
0}f_0^{\varepsilon}(x)=f_0(x),\quad
\lim\limits_{\varepsilon\rightarrow 0}u_0^{\varepsilon}(x)=u_0(x),$$
for almost all fixed $x\in\mathbb{R}^n;$
\item
for sufficiently small $\varepsilon$ and almost all fixed $(x,t)\in
{\mathbb R}^{n+1}$ every root $s_{k}$ of the equation
$\,u^\varepsilon_{0}(s)\,t+s-x=0\,$ belongs to the neighborhood
$U_{\bar s_{k}}(\varepsilon)$ of the root $\bar s_{k}$ of the
equation $\,\bar u^\varepsilon_{0}(s)\,t+s-x=0.$
\end{itemize}
\end{definition}

\begin{proposition}
The FP - solution $(f_{FP}(t,x),u_{FP}(t,x))$ to the Cauchy problem
(\ref{sist_pred1}), (\ref{sist_pred2}) with initial data
$(f_0(x),u_0(x))\in {\mathbb L}^2_{loc}({\mathbb R}^n)\cap {\mathbb
L}_{\infty}({\mathbb R}^n)$ does not depend of the choice of the
monotonic approximation $(f_0^{\varepsilon},u_0^{\varepsilon}).$
\end{proposition}

\textsc{Proof} Let us choose two the monotonic approximations
$\,(f_{01}^{\varepsilon},u_{01}^{\varepsilon})\,$ and
$(f_{02}^{\varepsilon},u_{02}^{\varepsilon})$ such that
$$\lim\limits_{\varepsilon\rightarrow
0}f^{\varepsilon}_{01}(x)=\lim\limits_{\varepsilon\rightarrow
0}f^{\varepsilon}_{02}(x)=f_0(x),$$
$$\lim\limits_{\varepsilon\rightarrow 0}u^{\varepsilon}_{01}(x)
=\lim\limits_{\varepsilon\rightarrow
0}u^{\varepsilon}_{02}(x)=u_0(x)$$ for almost all fixed
$x\in\mathbb{R}^n$. Then  the  couple
$$(f_0^{\varepsilon},u_0^{\varepsilon})=(f_{01}^{\varepsilon}-f_{02}^{\varepsilon},\,
u_{01}^{\varepsilon} -u_{02}^{\varepsilon})$$ can be considered as
initial data for the problem (\ref{sist_pred1})-(\ref{sist_pred2}).
To prove the proposition we have to show that the respective
solution  is identically zero almost everywhere.

Indeed, from (\ref{plotn}) we have for any $ t\ge 0,\,$ and almost
all $ x \in {\mathbb R}^n $
$$f_{FP}(t,x)=\sum\limits_{k}\,[\lim\limits_{\varepsilon\rightarrow
0}\left(\int\limits_{\mathbb{R}^n}f_{0}^{\varepsilon}(s)\delta(s-s^{\varepsilon}_{1,k}(t,x))ds\right)+$$
$$\lim\limits_{\varepsilon\rightarrow
0}\left(\int\limits_{\mathbb{R}^n}f_{02}^{\varepsilon}(s)\left(\delta(s-s^{\varepsilon}_{1,k}(t,x))\,-\,
\delta(s-s^{\varepsilon}_{2,k}(t,x))\right)ds\right)]=
$$
$$
\sum\limits_{k}\,[\lim\limits_{\varepsilon\rightarrow
0}(f^\varepsilon_{01}(s^\varepsilon_{1,k}(t,x))-f^\varepsilon_{02}(s^\varepsilon_{1,k}(t,x)))\,-\,
\lim\limits_{\varepsilon\rightarrow
0}(f^\varepsilon_{02}(s^\varepsilon_{2,k}(t,x))-f^\varepsilon_{02}(s^\varepsilon_{1,k}(t,x)))]
=0.$$ Here $s^\varepsilon_{i,k}(t,x)$ is the $k$-th solution
($k=1,2,...,K,\,\,i=1,2$) of equation
\begin{equation*}
u^\varepsilon_{0i}(s)\,t+s-x=0,\qquad i=1,2.
\end{equation*}
We have used the fact that
$|s^\varepsilon_{1,k}(t,x)\,-\,s^\varepsilon_{2,k}(t,x)|\,\to\,0,\,k=1,2,...,K$
as $\varepsilon\to 0.$
 Analogously proceeding from (\ref{sol_u_sdu}), we prove that
$u_{FP}(t,x)\equiv 0$ for almost all $ t\ge 0,\, x \in {\mathbb R}^n
$. $\square$

\bigskip

\section{The classical Riemann problem in the FP sense
for the 1D case}{\label{classRP}}

For the sake of simplicity we restrict ourselves to the
one-dimensional case and consider the following initial data:
\begin{equation}\label{f0_razrG}
f_0(x)=f_1(x)+\theta(x) f_2(x),
\end{equation}
\begin{equation}\label{u0_razrG}
u_0(x)=u_1(x)+\theta(x) u_2(x),
\end{equation}
where  $\theta$ is the Heaviside function with jump at $x_0$
(without loss of generality we shall assume $x_0=0$), $u_1,$ $u_2,$
$f_1,$ $f_2$ are continuously differentiable functions. We shall
dwell first on the case where $f_1, f_2, u_1, u_2$ are real
constants.
%\begin{equation}
%\label{f0_razr} f_0(x)=f_1+f_2\theta(x-x_0),
%\end{equation}
%\begin{equation}
%\label{u0_razr} u_0(x)=u_1+u_2\theta(x-x_0),
%\end{equation}
%where  $f_1, f_2, u_1, u_2$ are real constants.

According to Definition \ref{definit1} we must consider the smoothed
initial data instead of (\ref{f0_razrG}) and (\ref{u0_razrG}).  It
follows from Proposition \ref{prop1} that we can choose any couple
of smoothed initial data. It will be convenient to consider the
piecewise linear monotonic approximation of initial data of the form
\begin{equation}
\label{f} f^{\varepsilon}_0(x)=\begin{cases}
f_1,&\text{$x\leq -\varepsilon$,}\\
\dfrac{f_2}{2\varepsilon}x+f_1+\dfrac{f_2}2,&\text{$-\varepsilon<x<\varepsilon$,}\\
f_1+f_2,&\text{$x\geq \varepsilon$,}\\
\end{cases}
\end{equation}
\begin{equation}
\label{u} u^{\varepsilon}_0(x)=\begin{cases}
u_1,&\text{$x\leq -\varepsilon$,}\\
\dfrac{u_2}{2\varepsilon}x+u_1+\dfrac{u_2}2,&\text{$-\varepsilon<x<\varepsilon$,}\\
u_1+u_2,&\text{$x\geq \varepsilon$,}\\
\end{cases}
\end{equation}
where $f_1$, $f_2$, $u_1$ and $u_2$ are real constants.

From (\ref{plotn}) we can find the density
$\rho_\sigma^{\varepsilon}(t,x)$  corresponding to the smoothed
initial data $(f_0^{\varepsilon}(x),u_0^{\varepsilon}(x))$ (below we
shall omit the index $\sigma$ for short).

Let us set:
\begin{equation}
\label{E_hat1} \hat{E}_1(t,x,s)\equiv
\exp\left[-\dfrac{(s-x+u_1t)^2}{2\sigma^2t}\right],
\end{equation}

\begin{equation}
\label{E_hat2} \hat{E}_2(t,x,s)\equiv
\exp\left[-\dfrac{(s-x+(u_1+u_2)t)^2}{2\sigma^2t}\right],
\end{equation}

\begin{equation}
\label{E_hat3} \hat{E}_3(t,x,s)\equiv
\exp\left[-\dfrac{(s-x+(\dfrac{u_2}{2\varepsilon}s+u_1+\dfrac{u_2}2)t)^2}{2\sigma^2t}\right],
\end{equation}

Then

$$\rho^{\varepsilon}(t,x)=\dfrac1{\sqrt{2\pi t}\sigma}\left(\int\limits_{-\infty}^{-\varepsilon}f_1\hat{E}_1(t,x,s)d
s+\int\limits_{\varepsilon}^{+\infty}(f_1+f_2)\hat{E}_2(t,x,s)d
s+\right.$$

$$\left.+\int\limits_{-\varepsilon}^{\varepsilon}(\dfrac{f_2}{2\varepsilon}s+f_1+\dfrac{f_2}2)\hat{E}_3(t,x,s)d
s\right).$$

We denote the third integral  by $I_1^{\varepsilon,\sigma}$ and set
$p=\dfrac{(\dfrac{u_2t}{2\varepsilon}+1)s+(u_1+\dfrac{u_2}2)t-x}{\sigma\sqrt{t}}$.
We then have:

$$I_1^{\varepsilon,\sigma}=L^{\varepsilon}(t,x)\sigma\int\limits_{C_-^{\varepsilon}/\sigma\sqrt{t}}^{C_+^{\varepsilon}/\sigma\sqrt{t}}pe^{-\frac{p^2}2}d
p+F^{\varepsilon}(t,x)\sigma\int\limits_{C_-^{\varepsilon}/\sigma\sqrt{t}}^{C_+^{\varepsilon}/\sigma\sqrt{t}}e^{-\frac{p^2}2}d
p=$$

\begin{equation}
\label{I1}=L^{\varepsilon}(t,x)\sigma\left(e^{-\frac{(C_+^{\varepsilon})^2}{2t\sigma^2}}-e^{-\frac{(C_-^{\varepsilon})^2}{2t\sigma^2}}\right)+F^{\varepsilon}(t,x)\left[\Phi\left(\dfrac{C_+^{\varepsilon}}{\sigma}\right)-\Phi\left(\dfrac{C_-^{\varepsilon}}{\sigma}\right)\right],
\end{equation}

where
$\Phi(\alpha)=\dfrac1{\sqrt{2\pi}}\int\limits_{-\infty}^{\alpha}e^{-\frac{x^2}{2}}d
x$ is the Gauss function, $C_-^{\varepsilon}=u_1t-x-\varepsilon$,
$C_+^{\varepsilon}=(u_1+u_2)t-x+\varepsilon$, and

$$L^{\varepsilon}(t,x)=-\dfrac{\sqrt{2t}f_2\varepsilon}{\sqrt{\pi}(u_2t+2\varepsilon)^2},$$

\begin{equation}
\label{M_eps}F^{\varepsilon}(t,x)=\dfrac{2\varepsilon}{u_2t+2\varepsilon}\left(f_1+\dfrac{f_2}2+\dfrac{f_2}{u_2t+2\varepsilon}(x-(u_1+\dfrac{u_2}2)t)\right).
\end{equation}

It can easily be  seen that $\lim\limits_{\varepsilon\rightarrow
0}F^{\varepsilon}(t,x)=0$.

Finally, we get

\begin{equation}
\label{P1}\rho^{\varepsilon}(t,x)=f_1\Phi\left(\dfrac{C_-^{\varepsilon}}{\sigma\sqrt{t}}\right)+(f_1+f_2)\Phi\left(-\dfrac{C_+^{\varepsilon}}{\sigma\sqrt{t}}\right)+I_1^{\varepsilon,\sigma},
\end{equation}

 To find
$\hat{u}(t,x)$ we compute the numerator in formula
(\ref{sol_u_sdu}):
$$\dfrac1{\sqrt{2\pi t}\sigma}\int\limits_{\mathbb{R}}u_0^{\varepsilon}(s)f_0^{\varepsilon}(s)e^{-\frac{|u_0^{\varepsilon}(s)t+s-x|^2}{2\sigma^2t}}ds=u_1\rho^{\varepsilon}(t,x)+u_2(f_1+f_2)\Phi\left(-\dfrac{C_+^{\varepsilon}}{\sigma\sqrt{t}}\right)+I_2^{\varepsilon,\sigma},$$
where
\begin{equation}
\label{I2}I_2^{\varepsilon,\sigma}=G^{\varepsilon,\sigma}(t,x)+K^{\varepsilon}(t,x)
\sigma\left(e^{-\frac{(C_+^{\varepsilon})^2}{2t\sigma^2}}-e^{-\frac{(C_-^{\varepsilon})^2}{2t\sigma^2}}\right)+
\end{equation}
$$
N^{\varepsilon}(t,x)\left[\Phi\left(\dfrac{C_+^{\varepsilon}}{\sigma\sqrt{t}}\right)-\Phi\left(\dfrac{C_-^{\varepsilon}}{\sigma\sqrt{t}}\right)\right].
$$

Here we used the notation

$$K^{\varepsilon}(t,x)=\dfrac{u_2\sqrt{t}}{\sqrt{2\pi}(u_2t+2\varepsilon)}F^{\varepsilon}(t,x)+\dfrac{f_2\sqrt{t}}{2\sqrt{2\pi}\varepsilon}U^{\varepsilon}(t,x),$$

$$N^{\varepsilon}(t,x)=\left(\dfrac{u_2}{u_2t+2\varepsilon}(x-(u_1+\dfrac{u_2}2)t)+\dfrac{u_2}2\right)
F^{\varepsilon}(t,x).$$

We recall that $F^{\varepsilon}(t,x)$ was determined by
(\ref{M_eps}).

It is easy to deduce that $\lim\limits_{\varepsilon\rightarrow
0}N^{\varepsilon}(t,x)=0$.

Thus, we have the following result:

\begin{equation}
\label{sol_sdu1}\hat{u}^{\varepsilon}(t,x)=u_1
+\dfrac{u_2(f_1+f_2)\Phi\left(-\frac{C_+^{\varepsilon}}{\sigma\sqrt{t}}\right)+I_2^{\varepsilon,\sigma}}{f_1\Phi\left(\frac{C_-^{\varepsilon}}{\sigma\sqrt{t}}\right)+(f_1+f_2)\Phi\left(-\frac{C_+^{\varepsilon}}{\sigma\sqrt{t}}\right)+I_1^{\varepsilon,\sigma}},
\end{equation}
where $I_1^{\varepsilon,\sigma}$ and $I_2^{\varepsilon,\sigma}$ are
given by (\ref{I1}) and (\ref{I2}), respectively. Note that
$\dfrac{C_{\pm}^\varepsilon}{\sigma}\rightarrow$
$\rightarrow\pm\infty$ as $\sigma\rightarrow0$.

It can be checked that the initial conditions are satisfied, namely
$\rho(0,x)=f_0^{\varepsilon}(x)$ and
$\hat{u}^{\varepsilon}(0,x)=u_0^{\varepsilon}(x)$.

Now  we can find the generalized solution to the Riemann problem in
the form:
$$f_{FP}(t,x)=\lim\limits_{\varepsilon\rightarrow
0}(\lim\limits_{\sigma\rightarrow 0}\rho^{\varepsilon}(t,x)),$$
$$u_{FP}(t,x)=\lim\limits_{\varepsilon\rightarrow
0}(\lim\limits_{\sigma\rightarrow 0}\hat{u}^{\varepsilon}(t,x)).$$

Let us introduce the points
$\hat{x}^{\varepsilon}_1=u_1t-\varepsilon$ and
$\hat{x}^{\varepsilon}_2=(u_1+u_2)t+\varepsilon$. Their velocities
are $u_1$ and $u_1+u_2$, respectively.

\bigskip

We consider two cases:

{\bf 1.} $u_2>0$ (velocity of the point $\hat{x}^{\varepsilon}_2$ is
higher than the velocity of the point $\hat{x}^{\varepsilon}_1$). At
first, we can find
$f^{\varepsilon}(t,x)=\lim\limits_{\sigma\rightarrow
0}\rho^{\varepsilon}(t,x)$ from (\ref{P1}). It is easy to see that

$$
f^{\varepsilon}(t,x)=\begin{cases}
f_1,&\text{$x<\hat{x}^{\varepsilon}_1$,}\\
\dfrac{f_1}2-\dfrac12F^{\varepsilon}(t,x),&\text{$x=\hat{x}^{\varepsilon}_1$,}\\
F^{\varepsilon}(t,x),&\text{$\hat{x}^{\varepsilon}_1<x<\hat{x}^{\varepsilon}_2$,}\\
\dfrac{f_1+f_2}2+\dfrac12F^{\varepsilon}(t,x),&\text{$x=\hat{x}^{\varepsilon}_2$,}\\
f_1+f_2,&\text{$x>\hat{x}^{\varepsilon}_2$,}\\
\end{cases}
$$

 Let us note that this formula contains $F^{\varepsilon}(t,x)$ and
$\lim\limits_{\varepsilon\rightarrow 0}F^{\varepsilon}(t,x)=0.$
Thus, we obtain the following result for
$f_{FP}(t,x)=\lim\limits_{\varepsilon\rightarrow
0}f^{\varepsilon}(t,x)$:

$$
f_{FP}(t,x)=\begin{cases}
f_1,&\text{$x<u_1t$,}\\
\dfrac{f_1}2,&\text{$x=u_1t$,}\\
0,&\text{$u_1t<x<(u_1+u_2)t$,}\\
\dfrac{f_1+f_2}2,&\text{$x=(u_1+u_2)t$,}\\
f_1+f_2,&\text{$x>(u_1+u_2)t$.}\\
\end{cases}
$$
Further, from (\ref{sol_sdu1}) we find the solution of the gas
dynamic system with smooth initial data
$u^{\varepsilon}(t,x)=\lim\limits_{\sigma\rightarrow
0}\hat{u}_{\varepsilon}(t,x)$ as follows:
$$
u^{\varepsilon}(t,x)=\begin{cases}
u_1,&\text{$x<\hat{x}^{\varepsilon}_1$,}\\
u_1+\dfrac{N^{\varepsilon}(t,x)}{F^{\varepsilon}(t,x)},&\text{$\hat{x}^{\varepsilon}_1\leq x\leq \hat{x}^{\varepsilon}_2$,}\\
u_1+u_2,&\text{$x>\hat{x}^{\varepsilon}_2$.}\\
\end{cases}
$$
It can be shown that $$\lim\limits_{\varepsilon\rightarrow
0}\dfrac{N^{\varepsilon}(t,x)}{F^{\varepsilon}(t,x)}=\lim\limits_{\varepsilon\rightarrow
0}\dfrac{u_2}{u_2t+2\varepsilon}(x-(u_1+\dfrac{u_2}2)t)+\dfrac{u_2}2=\dfrac{x}{t}-u_1.$$

\begin{figure}[h]
\begin{minipage}{0.5\columnwidth}
\centerline{\includegraphics[width=1\columnwidth]{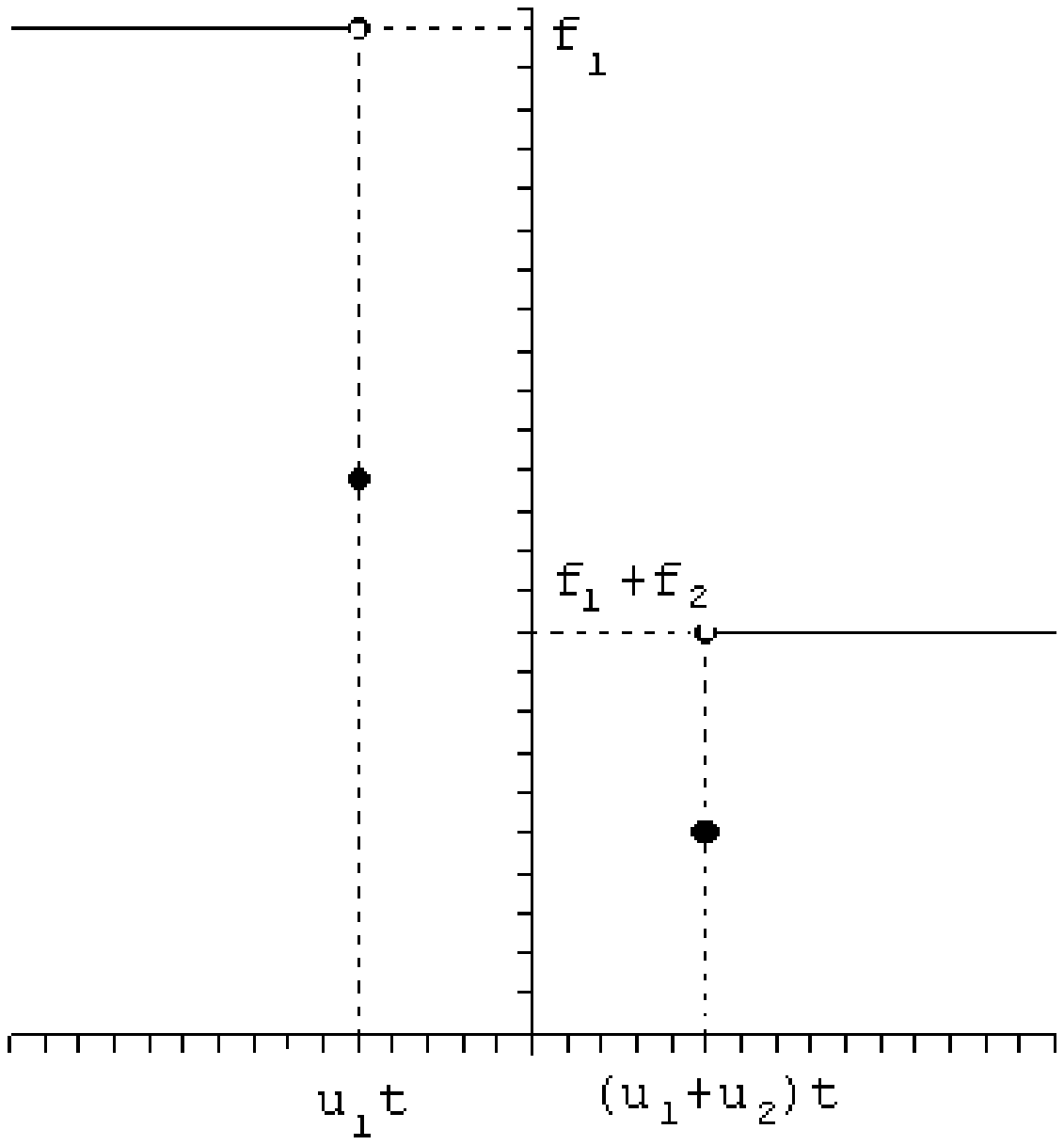}}
\end{minipage}%
%\end{figure}%
%\begin{figure}[h]
\begin{minipage}{0.5\columnwidth}
\centerline{\includegraphics[width=1\columnwidth]{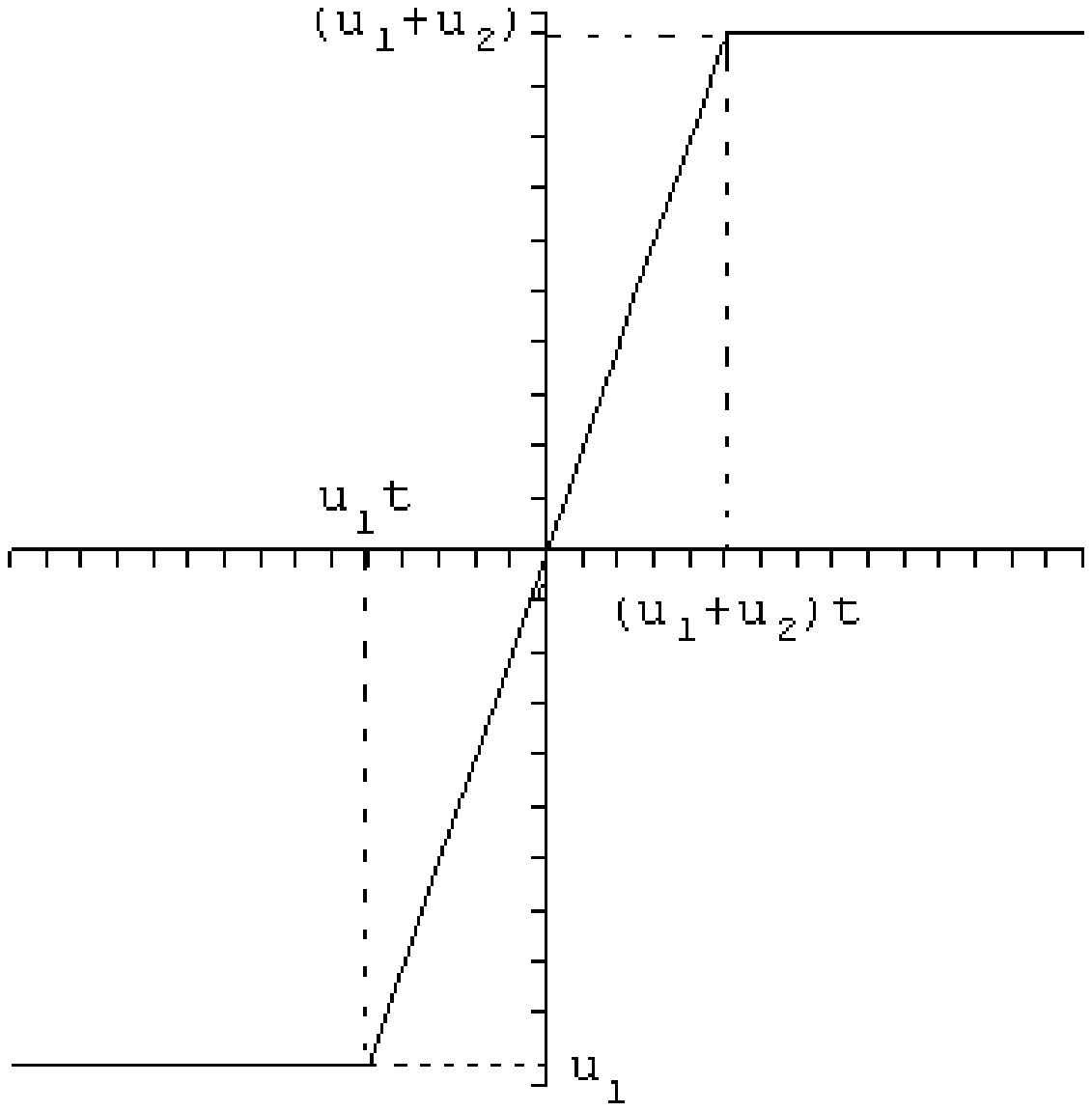}}
\end{minipage}
\caption{Density and velocity, $u_2>0.$} \label{Fig1}
\end{figure}

Thus, we get the following solution for
$u(t,x)=\lim\limits_{\varepsilon\rightarrow
0}u^{\varepsilon}(t,x)$:
$$
u_{FP}(t,x)=\begin{cases}
u_1,&\text{$x< u_1t$,}\\
\dfrac{x}{t},&\text{$u_1t\leq x\leq(u_1+u_2)t$,}\\
u_1+u_2,&\text{$x>(u_1+u_2)t$.}\\
\end{cases}
$$
We can see that the velocity includes the rarefaction wave (see
Fig.\ref{Fig1}). This is a well known self-similar solution to the
Riemann problem with constant left-hand and right-hand states for
the Burgers equation (\cite{ShengZhang}).

% According to Proposition \ref{prop1} and Remark \ref{remf0_jump}, the pair $(f,u)$ solve the
% pressureless gas dynamics system almost everywhere on ${\mathbb R}_+\times {\mathbb R}^n.$

It is interesting to note that if we first compute the limit in
$\varepsilon$  we get the solution
$u(t,x)=u_1+u_2\theta(x-(u_1+\dfrac{u_2}2)t)$, which is unstable
with respect to  small perturbations.

\begin{figure}[h]
\begin{minipage}{0.5\columnwidth}
\centerline{\includegraphics[width=1\columnwidth]{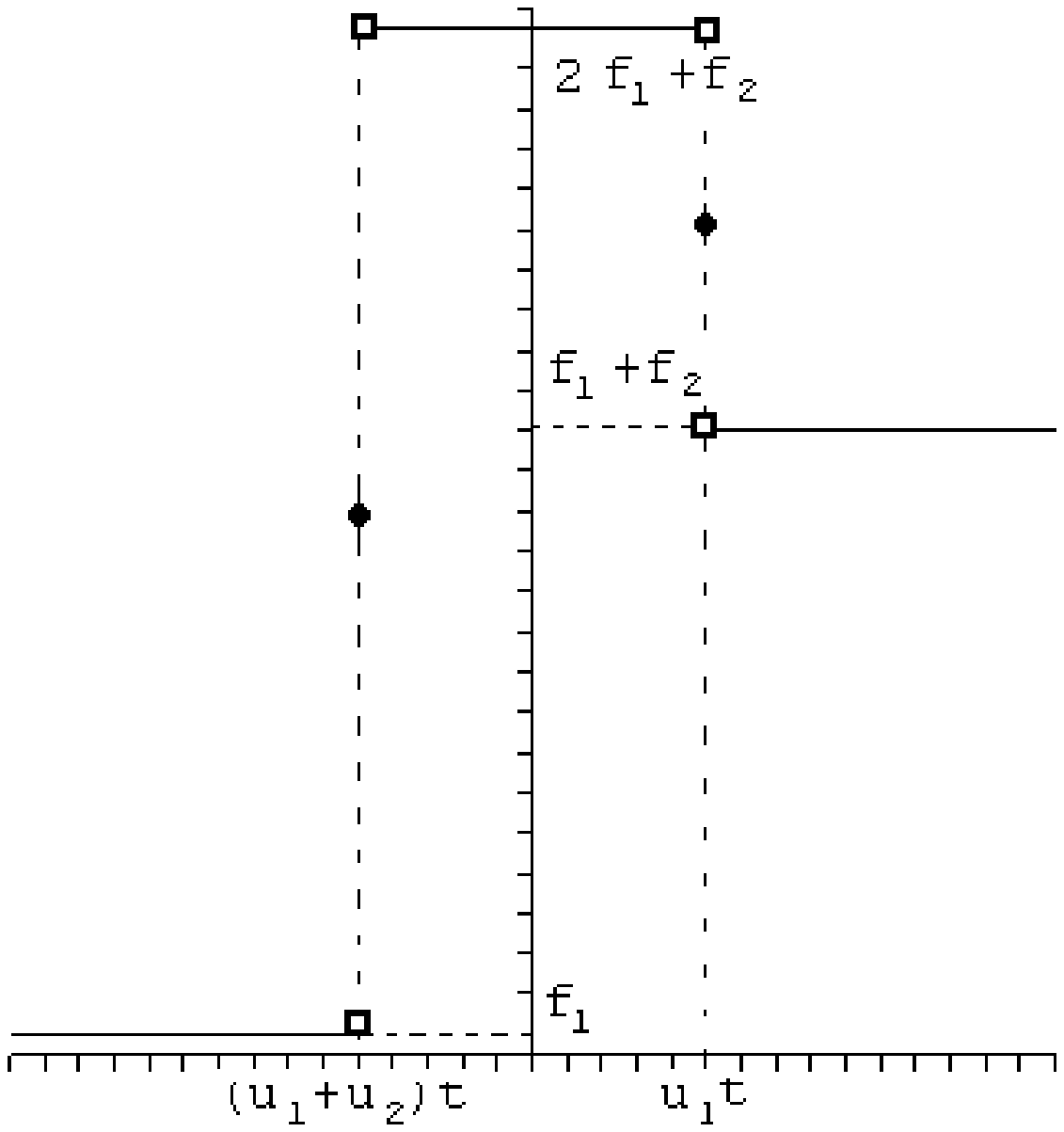}}
\end{minipage}%
%\end{figure}%
%\begin{figure}[h]
\begin{minipage}{0.5\columnwidth}
\centerline{\includegraphics[width=1\columnwidth]{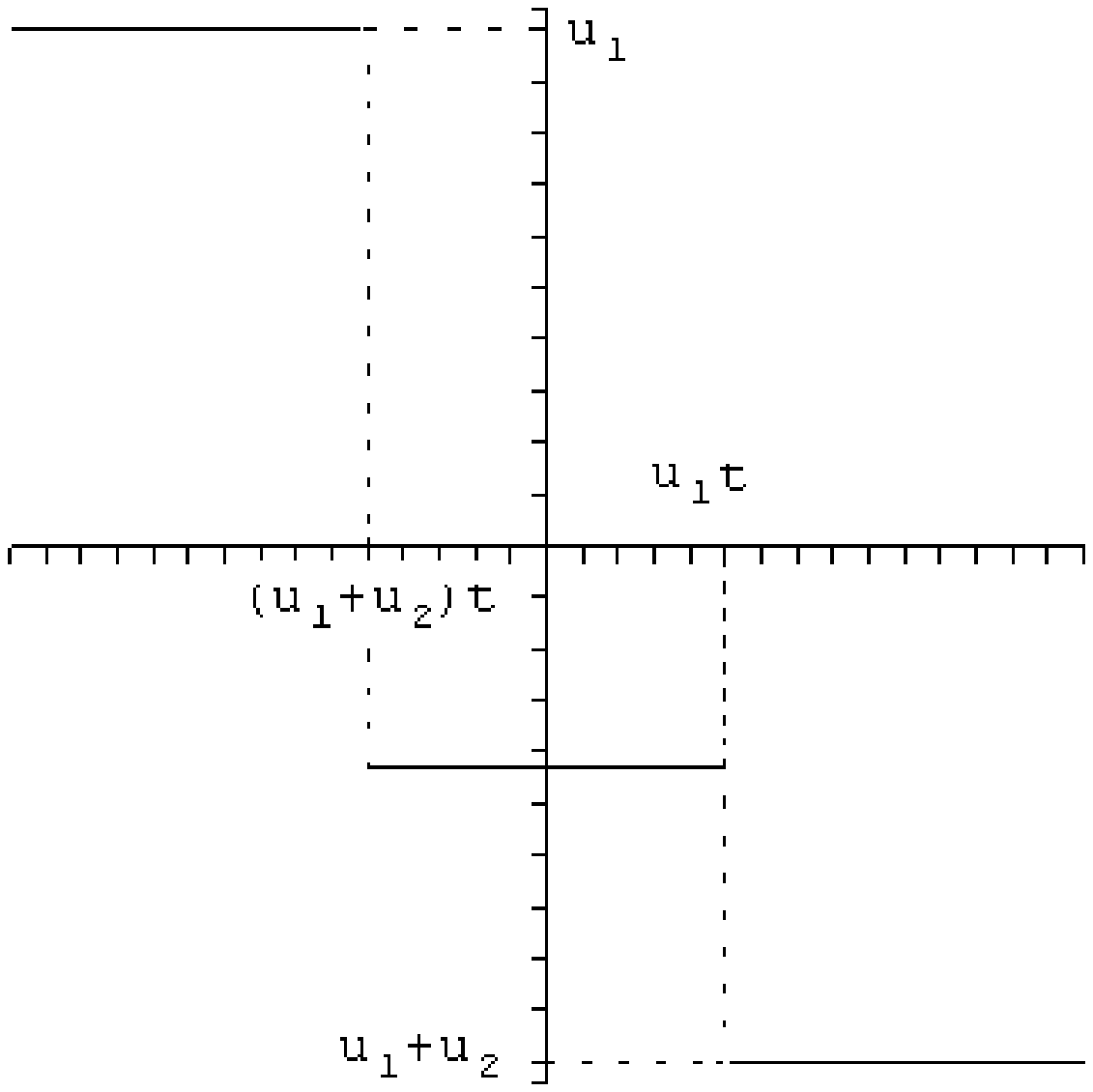}}
\end{minipage}
\caption{Density and velocity, $u_2<0.$} \label{Fig2}
\end{figure}

{\bf 2.} $u_2<0$ (the velocity of $\hat{x}^{\varepsilon}_2$ is
higher than the velocity of $\hat{x}^{\varepsilon}_1$). From
(\ref{P1}) and (\ref{sol_sdu1}) we find as before (see
Fig.\ref{Fig2})
$$
f_{FP}(t,x)=\begin{cases}
f_1,&\text{$x<(u_1+u_2)t$,}\\
\dfrac{3f_1+f_2}2,&\text{$x=(u_1+u_2)t$,}\\
2f_1+f_2,&\text{$(u_1+u_2)t<x<u_1t$,}\\
\dfrac{3f_1+2f_2}2,&\text{$x=u_1t$,}\\
f_1+f_2,&\text{$x>u_1t$,}\\
\end{cases}
$$
$$
u_{FP}(t,x)=\begin{cases}
u_1,&\text{$x<(u_1+u_2)t$,}\\
u_1+\dfrac{f_1+f_2}{2f_1+f_2}u_2,&\text{$(u_1+u_2)t\leq x\leq u_1t$,}\\
u_1+u_2,&\text{$x> u_1t$.}\\
\end{cases}
$$

%\begin{remark}\,
% We can extend our results to the case of any initial distribution
%$f_0(x)$, e.g. a Gaussian distribution.
%\end{remark}

\section{The singular Riemann problem in the FP sense}{\label{singRP}

By the singular Riemann problem we mean the Cauchy problem with the
following data:
\begin{equation}
\label{f0_razrS}
f_0(x)=f_1+f_2\theta(x-x_0)+f_3\delta(x-x_0)=f_{reg}^0+f_{sing}^0,
\end{equation}
\begin{equation}
\label{u0_razrS} u_0(x)=u_1+u_2\theta(x-x_0),
\end{equation}
which differs from the classical Riemann problem (\ref{f0_razrG}),
(\ref{u0_razrG}) by a singular part $f_{sing}^0,$ the $\delta$ -
function of constant amplitude $f_3$ concentrated at the jump of the
density. We again set $x_0$ equal to zero (without loss of
generality).

%As we will show in  Sec.\ref{Arising}, the $\delta$ - function in
%the density component  is intrinsic to the formation of a
%singularity from smooth initial data, therefore for the pressureless
%gas dynamics system the singular Riemann problem is more natural
%than the classical  one.

Thus, we should adapt the definition of the generalized solution to
the case of a singular density. As before we want to approximate
initial data by  smooth functions,  the $\delta$-function will be
instead naturally approximated in the space $\mathcal D'(\mathbb R)$
of distributions.

\begin{definition}\label{d1} We call the couple of functions $(f_{FP}(t,x),u_{FP}(t,x))$
a generalized solution to the problem (\ref{sist_pred1}) ,
(\ref{sist_pred2}), (\ref{f0_razrS}), (\ref{u0_razrS}) in the sense
of free particles (FP), if
$$f_{FP}(t,x)=\lim\limits_{\varepsilon\rightarrow
0}(\lim\limits_{\sigma\rightarrow
0}\rho_{\sigma,reg}^{\varepsilon}(t,x))\,+\,\lim\limits_{\varepsilon\rightarrow
0}(\lim\limits_{\sigma\rightarrow
0}\rho_{\sigma,sing}^{\varepsilon}(t,x))\,=\,f_{reg}(t,x)\,+\,f_{sing}(t,x),$$
$$u_{FP}(t,x)=\lim\limits_{\varepsilon\rightarrow
0}(\lim\limits_{\sigma\rightarrow
0}\hat{u}_\sigma^{\varepsilon}(t,x)),$$ where the limits are meant
as pointwise at almost all points except for the singular part
$$\lim\limits_{\varepsilon\rightarrow
0}(\lim\limits_{\sigma\rightarrow 0}\rho_{\sigma,
sing}^{\varepsilon}(t,x)),$$ where the outer limit in $\varepsilon$
is meant in the sense of distributions.

Here $(\rho_{\sigma, reg}^{\varepsilon}(t,x)\,+\,\rho_{\sigma,
sing}^{\varepsilon}(t,x),\hat{u}^{\varepsilon}(t,x))$ satisfy the
system (\ref{sist_obw1D1}), (\ref{sist_obw1D2}) with initial data
$(f_{0,reg}^{\varepsilon}(x)\,+\,f_{0,sing}^{\varepsilon}(x),u_0^{\varepsilon}(x))$
from the class $C^1(\mathbb{R}^n)$ such that
$$\lim\limits_{\varepsilon\rightarrow
0}f_{0,reg}^{\varepsilon}(x)=f_0^{reg}(x),\quad
\lim\limits_{\varepsilon\rightarrow 0}u_0^{\varepsilon}(x)=u_0(x)$$
for almost all fixed $x\in\mathbb{R}^n$ and
$$\lim\limits_{\varepsilon\rightarrow
0}f_{0,sing}^{\varepsilon}(x)=f_0^{sing}(x),$$ the limit in
$\varepsilon$ being in the sense of distributions.
\end{definition}

\bigskip

We are going to solve the singular Riemann problem with constant
left and right states.  To approximate the $\delta$-function we use
the well known fact that
$$\dfrac{1}{\sqrt{2\pi\varepsilon}}\,\,e^{-\frac{x^2}{2\varepsilon}}\rightarrow\delta(x),\quad \varepsilon\rightarrow 0\quad \text{in}\,\,
\mathcal D'(\mathbb R).$$

The part of the solution that relates to the regular part of the
density ($f_3=0$) is found in Sec.\ref{classRP}. Now we have to
calculate the singular part of the density
$\rho^{\varepsilon}_{\sigma, sing}(t,x)$ for
$f_{0,sing}^{\varepsilon}(x)=\,f_3\,\dfrac1{\sqrt{2\pi\varepsilon}}e^{-\frac{x^2}{2\varepsilon}}$.
As before we omit the index $\sigma.$
%After this we will add it to the previous result.

From (\ref{plotn}) we obtain:

$$\rho^{\varepsilon}_{sing}(t,x)=\dfrac{f_3}{\sqrt{2\pi(\sigma^2t+\varepsilon)}}\left(\exp\left[-\dfrac{(u_1t-x)^2}{2(\sigma^2t+\varepsilon)}\right]\Phi\left(\dfrac{D_1^{\sigma,\varepsilon}}{\sigma\sqrt{t}}\right)+\right.$$
\begin{equation}\label{rho_sing}\left.+\exp\left[-\dfrac{((u_1+u_2)t-x)^2}{2(\sigma^2t+\varepsilon)}
\right]\Phi\left(-\dfrac{D_2^{\sigma,\varepsilon}}{\sigma\sqrt{t}}\right)\right)+J_2^{\varepsilon,\sigma}=\end{equation}
$$=\dfrac{f_3}{\sqrt{2\pi(\sigma^2t+\varepsilon)}}\left(G_1^{\sigma,\varepsilon}(t,x)\Phi\left(\dfrac{D_1^{\sigma,\varepsilon}}{\sigma\sqrt{t}}\right)+G_2^{\sigma,\varepsilon}(t,x)\Phi\left(-\dfrac{D_2^{\sigma,\varepsilon}}{\sigma\sqrt{t}}\right)\right)+J_2^{\varepsilon,\sigma},$$
where

$$D_1^{\sigma,\varepsilon}=\dfrac{(u_1t-x)\sqrt{\varepsilon}}{\sqrt{\sigma^2t+\varepsilon}}-\sqrt{\varepsilon(\sigma^2t+\varepsilon)},$$

$$D_2^{\sigma,\varepsilon}=\dfrac{((u_1+u_2)t-x)\sqrt{\varepsilon}}{\sqrt{\sigma^2t+\varepsilon}}+\sqrt{\varepsilon(\sigma^2t+\varepsilon)},$$
and
$J_2^{\varepsilon,\sigma}=\int\limits_{-\varepsilon}^{\varepsilon}\exp\left[-\dfrac{s^2}{2\varepsilon}-\dfrac{(s-x+(\dfrac{u_2}{2\varepsilon}s+u_1+\dfrac{u_2}2)t)^2}{2\sigma^2t}\right]ds$.

It can be checked that

$$J_2^{\varepsilon,\sigma}=\dfrac{f_3}{\sqrt{2\pi}\sqrt{A^{\sigma,\varepsilon}}}\exp\left[-\dfrac{((u_1+\frac{u_2}2)t-x)^2}{2A^{\sigma,\varepsilon}}\right]\left(\Phi\left(\dfrac{K_+^{\sigma,\varepsilon}}{\sigma\sqrt{t}}\right)-\Phi\left(\dfrac{K_-^{\sigma,\varepsilon}}{\sigma\sqrt{t}}\right)\right),$$

where

\begin{equation}
\label{A_nu}A^{\sigma,\varepsilon}=\sigma^2t+\varepsilon(\frac{u_2t}{2\varepsilon}+1)^2,
\end{equation}

\begin{equation}
\label{B_nu}B^{\sigma,\varepsilon}=\dfrac{\varepsilon(\frac{u_2t}{2\varepsilon}+1)((u_1+\frac{u_2}2)t-x)}{\sqrt{A^{\sigma,\varepsilon}}},
\end{equation}

\begin{equation}
\label{K_nu}K_{\pm}^{\sigma,\varepsilon}=\dfrac1{\sqrt{\varepsilon}}\left(\sqrt{A^{\sigma,\varepsilon}}(\pm\varepsilon)+B^{\sigma,\varepsilon}\right).
\end{equation}

It is easy to see that
$$D_1^{\varepsilon}=\lim\limits_{\sigma\rightarrow
0}D_1^{\sigma,\varepsilon}=u_1t-x-\varepsilon,\quad
D_2^{\varepsilon}=\lim\limits_{\sigma\rightarrow
0}D_2^{\sigma,\varepsilon}=(u_1+u_2)t-x-\varepsilon,$$$$
K_{\pm}^{\varepsilon}=\lim\limits_{\sigma\rightarrow
0}K_{\pm}^{\sigma,\varepsilon}(\varepsilon)=(\dfrac{u_2t}{2\varepsilon}+1)(\pm\varepsilon)+(u_1+\dfrac{u_2}2)t-x,$$
i.e. $K_-^{\varepsilon}=D_1^{\varepsilon}=C_-^{\varepsilon}$,
$K_+^{\varepsilon}=D_2^{\varepsilon}=C_+^{\varepsilon}$. Therefore,
in this case we also have two jump points
$\hat{x}_1^\varepsilon=u_1t-\varepsilon$ and
$\hat{x}_2^\varepsilon=(u_1+u_2)t+\varepsilon$.

We have to consider two cases $u_2>0$ and $u_2<0,$ as before. But
for $f_{0}^{sing}(x)=$ $=f_3\,\delta(x)$ we obtain the same result:
$$f_{sing}(t,x)=\lim\limits_{\varepsilon\rightarrow
0}\lim\limits_{\sigma\rightarrow
0}\rho^{\varepsilon}_{sing}(t,x)=\dfrac12 f_3(\delta(x
-u_1t)+\delta(x-(u_1+u_2)t)).$$  Thus, we have two
$\delta$-functions with equal amplitudes $\dfrac12F(t)$, which move
with the jump points.

Further,  we can find the regular part of the density
$f_{reg}(t,x)=\lim\limits_{\varepsilon\rightarrow
0}\lim\limits_{\sigma\rightarrow 0}\rho_{reg}^{\varepsilon}(t,x)$,
where
\begin{equation}
\label{P_delta}\rho_{reg}^{\varepsilon}(t,x)=f_1\Phi\left(\dfrac{C_-^{\varepsilon}}
{\sigma\sqrt{t}}\right)+(f_1+f_2)\Phi\left(-\dfrac{C_+^{\varepsilon}}{\sigma\sqrt{t}}\right)+I_1.
%\rho^{\varepsilon}_{\delta}(t,x).
\end{equation}

Analogously, from (\ref{sol_u_sdu}) we can calculate the velocity
$\hat{u}(t,x)$:

$$\hat{u}^{\varepsilon}(t,x)=u_1+\dfrac{u_2(f_1+f_2)\Phi\left(-\frac{C_+^{\varepsilon}}{\sigma\sqrt{t}}\right)+I_2^{\varepsilon,\sigma}+\dfrac{f_3u_2}{\sqrt{2\pi(\sigma^2t+\varepsilon)}}G_2^{\sigma,\varepsilon}(t,x)\Phi\left(-\frac{D_2^{\sigma,\varepsilon}}{\sigma\sqrt{t}}\right)+J_3^{\varepsilon,\sigma}}{\rho^{\varepsilon}(t,x)},$$

where

$$J_3^{\varepsilon,\sigma}=\dfrac{f_3}{\sqrt{2\pi}}G^{\sigma,\varepsilon}(t,x)\left(\dfrac{u_2\sqrt{t\varepsilon}}{2\varepsilon
A^{\sigma,\varepsilon}}\sigma\left[e^{-\frac{(K^{\sigma,\varepsilon}_+)^2}{2t\sigma^2}}-e^{-\frac{(K^{\sigma,\varepsilon}_-)^2}{2t\sigma^2}}\right]+\right.$$

$$\left.+\dfrac1{\sqrt{A^{\sigma,\varepsilon}}}\left(\dfrac{u_2}2-\dfrac{u_2B^{\sigma,\varepsilon}}{2\varepsilon A^{\sigma,\varepsilon}}\right)\left[\Phi\left(\dfrac{K_+^{\sigma,\varepsilon}}{\sigma\sqrt{t}}\right)-\Phi\left(\dfrac{K_-^{\sigma,\varepsilon}}{\sigma\sqrt{t}}\right)\right]\right),$$

$$G^{\sigma,\varepsilon}(t,x)=\exp\left[-\dfrac{((u_1+\frac{u_2}2)t-x)^2}{2A^{\sigma,\varepsilon}}\right],$$

$$G_2^{\sigma,\varepsilon}(t,x)=\exp\left[-\dfrac{((u_1+u_2)t-x)^2}{2(\sigma^2t+\varepsilon)}\right].$$

The expressions for
$\rho^{\varepsilon}(t,x)\,=\,\rho_{reg}^{\varepsilon}(t,x)\,+\,\rho_{sing}^{\varepsilon}(t,x)$
and $I_2^{\varepsilon,\sigma}$  have been given in (\ref{P_delta}),
(\ref{rho_sing}) and (\ref{I2}), respectively.

\textbf{1.} $u_2>0$ then:

$$ u^{\varepsilon}(t,x)=\begin{cases}
u_1,&\text{$x<\hat{x}^{\varepsilon}_1$,}\\
u_1+\dfrac{\dfrac{f_3}{\sqrt{2\pi}}T^{\varepsilon}G^{\varepsilon}(t,x)+N^{\varepsilon}}{\dfrac{f_3}{\sqrt{2\pi\varepsilon}}\dfrac{2\varepsilon}{2\varepsilon+u_2t}G^{\varepsilon}(t,x)+F^{\varepsilon}},&\text{$\hat{x}^{\varepsilon}_1\leq x\leq\hat{x}^{\varepsilon}_2$,}\\
u_1+u_2,&\text{$x>\hat{x}^{\varepsilon}_2$,}\\
\end{cases}
$$

where

$$T^{\varepsilon}=\lim\limits_{\sigma\rightarrow
0}\dfrac1{\sqrt{A^{\sigma,\varepsilon}}}\left(\dfrac{u_2}2-\dfrac{u_2B^{\sigma,\varepsilon}}{2\varepsilon
A^{\sigma,\varepsilon}}\right)=\dfrac{u_2}2-\dfrac{2\varepsilon((u_1+\dfrac{u_2}2)t-x)u_2}{(2\varepsilon+u_2t)^2\sqrt{\varepsilon}},$$

Therefore

$$ u_{FP}(t,x)=\lim\limits_{\varepsilon\rightarrow
0}u^{\varepsilon}(t,x)=\begin{cases}
u_1,&\text{$x<u_1t$,}\\
\dfrac{x}{t},&\text{$u_1t\leq x\leq(u_1+u_2)t$,}\\
u_1+u_2,&\text{$x>\hat{x}^(u_1+u_2)t$.}\\
\end{cases}
$$
Further, from  (\ref{P_delta}), (\ref{rho_sing}) we find  the
density as follows:
$$
f_{FP}(t,x)=\dfrac12f_3(\delta(x-u_1t)+\delta(x-(u_1+u_2)t))+\begin{cases}
f_1,&\text{$x<u_1t$,}\\
\dfrac{f_1}2,&\text{$x=u_1t$,}\\
0,&\text{$u_1t<x<(u_1+u_2)t$,}\\
\dfrac{f_1+f_2}2,&\text{$x=(u_1+u_2)t$,}\\
f_1+f_2,&\text{$x>(u_1+u_2)t$,}\\
\end{cases}
$$

\textbf{2.} If $u_2<0,$ we have

$$ u^{\varepsilon}(t,x)=\begin{cases}
u_1,&\text{$x<\hat{x}^{\varepsilon}_1$,}\\
u_1+\dfrac{u_2(f_1+f_2)-N^{\varepsilon}+\dfrac{f_3}{\sqrt{2\pi\varepsilon}}(G_2^{\varepsilon}-\varepsilon T^{\varepsilon}G^{\varepsilon})}{2f_1+f_2+F^{\varepsilon}+\dfrac{f_3}{\sqrt{2\pi\varepsilon}}(G_1^{\varepsilon}+G_2^{\varepsilon}-\dfrac{2\varepsilon}{2\varepsilon+u_2t}G^{\varepsilon})},&\text{$\hat{x}^{\varepsilon}_1\leq x\leq\hat{x}^{\varepsilon}_2$,}\\
u_1+u_2,&\text{$x>\hat{x}^{\varepsilon}_2$,}\\
\end{cases}
$$

where $G_1^{\varepsilon}\equiv
\exp\left[-\dfrac{(u_1t-x)^2}{2\varepsilon}\right]$,
$G_2^{\varepsilon}\equiv
\exp\left[-\dfrac{((u_1+u_2)t-x)^2}{2\varepsilon}\right]$.

Thus,

$$ u_{FP}(t,x)=\begin{cases}
u_1,&\text{$x<(u_1+u_2)t$,}\\
u_1+\dfrac{(f_1+f_2)}{2f_1+f_2}u_2,&\text{$(u_1+u_2)t\leq x\leq u_1t$,}\\
u_1+u_2,&\text{$x>u_1t$.}\\
\end{cases}
$$

As a consequence, for the velocity we obtain the same result as in
the non-singular case. Analogously,
$$
f_{FP}(t,x)=\dfrac12f_3(\delta(x-u_1t)+\delta(x-(u_1+u_2)t))+\begin{cases}
f_1,&\text{$x<(u_1+u_2)t$,}\\
\dfrac{3f_1+f_2}2,&\text{$x=(u_1+u_2)t$,}\\
2f_1+f_2,&\text{$(u_1+u_2)t<x<u_1t$,}\\
\dfrac{3f_1+2f_2}2,&\text{$x=u_1t$,}\\
f_1+f_2,&\text{$x>u_1t$.}\\
\end{cases}
$$

\section{Singularity arising from smooth
data}\label{Arising}

We are going to show that at the point of formation of a singularity
from smooth initial data in the solution to the pressureless gas
dynamics model a $\delta$ -- function appears in the density
component.  For the sake of simplicity we again restrict  ourselves
to the 1D case.
\begin{Theorem}\label{T1}
{\rm (Asymptotics of the approximating solution)} Let the initial
data $\,(f_0,\,u_0)\, $ for the system (\ref{sist_pred1}),
(\ref{sist_pred2})  be at least $C^m$ -- smooth and bounded, $m\ge
2.$ Assume that there exists an instant $0< t_*<\infty,$ such that
$t_*=\inf\limits_{x\in \mathbb R}\left(-\frac{1}{u'_0(x)}\right) $
and $u^{(k)}(s_*)=0, \,k=1,...,m-1,$ however $u^{(m)}(s_*)$ does not
vanish at the point $s_*(t_*,x_*),$ where $s_*$ is a solution to the
equation $u_0(s)\,t_*\,=\,x_*-\,s.$ Here $\,x_* $ is such that the
line $y=\frac{x_*-s}{t_*}$ intersects the graph of the initial
velocity $y=u_0(s)$ at a unique point and it is tangent to the
graph.

Then at the moment $t=t_*$   the following properties of the
function $\rho_\sigma (t,x),$ entering the solution $(\rho_\sigma,
\hat u_\sigma)$  to the system (\ref{sist_obw1D1}),
(\ref{sist_obw1D2}) hold:
\begin{itemize}
\item
\begin{equation}\label{T1i}
\rho_\sigma(t_*,x_*)\,\sim\,B(x_*,t_*)\,f_0(s_*)\,\sigma^{-\frac{m-1}{m}},\,\quad\sigma\to
0,
\end{equation}
where
$$
B(x_*,t_*)\,=\,K_m\,
\frac{|u_0'(s_*)|^{\frac{m-1}{2m}}}{|u_0^{(m)}(s_*)|^{\frac{1}{m}}},\qquad
K_m\,=\,\frac{1}{2^{\frac{m-1}{2m}}\,m\,\sqrt{\pi}}\,(m!)^{\frac{1}{m}}\,\Gamma\left(\frac{1}{2m}\right),
$$
$\sim$ means that the quotient of the left-hand side by the
right-hand side converges to $1$ as $\sigma\to 0;$
\item
\begin{equation}\label{T1ii}
\rho_\sigma(t_*,x)\,\to \,f_0(s_0(t_*,x)),\quad\sigma\to 0, \quad
\mbox{for}\,
 x\ne x_*,
\end{equation}
where the function $s_0(t,x)$ has been introduced in the proof of
Proposition \ref{prop1}.
\end{itemize}
\end{Theorem}

\proof  Proceeding as in the proof of Proposition \ref{prop1} we can
readily obtain the property (\ref{T1ii}). Thus, let us dwell on the
property  (\ref{T1i}).

Let us analyze the formula (\ref{plotn}) at the point $(t_*, x_*)$.
To that end we note that since $\,u_0(s_*)\,t_*=x_*-s_*\,$ and
$\,t_*=-\frac{1}{u_0(s_*)},\,$ for $\,s\,$ belonging to the
$\varepsilon$ - neighborhood $U_{s_*}(\varepsilon)$ of the point
$s_*$ we have
$$
u_0(s)\,t_*+s-x\,=\,(u_0(s_*)+u'_0(s_*)\,(s-s_*)+\frac{1}{m!}\,u^{(m)}_0(s_{**)}\,(s-s_*)^m)\,t_*\,+s-x=$$
\begin{equation}\label{U_eps}
=\frac{1}{m!}\,u^{(m)}_0(s_{**)}\,(s-s_*)^m\,t_*-(x-x_*),
\end{equation}
with $s_{**}\in U_{s_*}(\varepsilon),$ $\,\varepsilon>0.$

Then from (\ref{plotn})  we have
$$\rho (t_*,x_*)\,=\,
\dfrac1{\sqrt{2\pi t_*}\sigma}\,\int\limits_{\mathbb
R}\,f_0(s)e^{-\frac{(u_0(s)t+s-x_*)^2}{2\sigma^2t_*}}\,ds\,=$$
$$
\dfrac1{\sqrt{2\pi t_*}\sigma}\,\int\limits_{\mathbb
R}\,f_0(s)\,e^{-\frac{(\frac{u_0^{(m)}(s_*)t_*}{m!}\,(s-s_*)^m)^2}{2\sigma^2t_*}}\,ds\,+
$$
$$
\dfrac1{\sqrt{2\pi t_*}\sigma}\,\int\limits_{\mathbb R}\,f_0(s)
\left(e^{-\frac{(u_0(s)t+s-x_*)^2}{2\sigma^2t_*}}\,-\,
e^{-\frac{(\frac{u_0^{(m)}(s_*)t_*}{m!}\,(s-s_*)^m)^2}{2\sigma^2t_*}}\right)\,
ds=$$
$$I_1+I_2.$$
The first integral $I_1$ is equal to
$$\sigma^{-\frac{m-1}{m}}\,I_{\sigma}(t_*,x_*),$$
where $I_{\sigma}(t_*,x_*)\to B(x_*,t_*)\,f_0(s_*) $ as $\sigma\to
0,$ where $B(x_*,t_*)$ is specified in the statement of Theorem
\ref{T1}. To evaluate $I_\sigma$ we have used the formula
$$
\int\limits_{{\mathbb
R}}\,e^{-a^2\,s^{2m}}\,ds\,=\,\frac{\Gamma\left(\frac{1}{2m}\right)}{|a|^{\frac{1}{m}}\,m},\quad
a={\rm const}\ne 0.
$$
Now let as prove that $I_2\to 0$ as $\sigma\to 0.$

Let us choose $\varepsilon >0$ so small that for all $s\in
U_{s_*}(\varepsilon)$
$$
\big|\,e^{-\frac{(u_0(s)t+s-x_*)^2}{2\sigma^2t_*}}\,-\,
e^{-\frac{(\frac{u_0^{(m)}(s_*)t_*}{m!}\,(s-s_*)^m)^2}{2\sigma^2t_*}}\,\big|<
\sigma.
$$

Then
$$
|I_2|\le \dfrac1{\sqrt{2\pi
t_*}\sigma}\,\int\limits_{U_{s_*}(\varepsilon)}\,|f_0(s)|
\left|e^{-\frac{(u_0(s)t+s-x_*)^2}{2\sigma^2t_*}}\,-\,
e^{-\frac{(\frac{u_0^{(m)}(s_*)t_*}{m!}\,(s-s_*)^m)^2}{2\sigma^2t_*}}\right|\,
ds\,+$$
$$
\dfrac1{\sqrt{2\pi t_*}\sigma}\,\int\limits_{{\mathbb R}\backslash
U_{s_*}(\varepsilon)}\,|f_0(s)|
\left(e^{-\frac{(u_0(s)t+s-x_*)^2}{2\sigma^2t_*}}\,+\,
e^{-\frac{(\frac{u_0^{(m)}(s_*)t_*}{m!}\,(s-s_*)^m)^2}{2\sigma^2t_*}}\right)\,
ds.
$$
The first part in the right-hand side of the inequality due to the
boundedness of $f_0$ is less than $\,{\rm const}\cdot\varepsilon,\,
$ the second part  tends to zero as $\sigma\to 0$ due to the
boundness of $f_0.$ Since $\varepsilon $ can be chosen arbitrarily
small, the statement is proved. $\square$

\begin{remark}
The following asymptotics of $K_m$  holds:
$$
K_m\,=\,\frac{\sqrt{2}}{e}\,m,+\,O(\ln m)\quad m\to\infty.
$$
\end{remark}

\begin{Theorem} {\rm(Amplitude of the $\delta-$ function)}\label{T2}
Let the initial data $\,(f_0,\,u_0)\, $ for the system
(\ref{sist_pred1}), (\ref{sist_pred2})  be  $C^1$ -- smooth and
bounded and let the critical instant $t_*=\inf\limits_{x\in \mathbb
R}\left(-\frac{1}{u'_0(x)}\right) $ be positive and finite. Assume
that the initial datum $u_0$ is linear on the segment
$\Omega=(x_1,x_2),$ moreover, the second left-hand derivative
$u_0''(x_1-0)$ at the point $x_1$ and the second right-hand
derivative $u_0''(x_2+0)$ at the point $x_2$ do not vanish. Let
$x_*$ be the unique point such that the line  $y=\frac{x_*-s}{t_*}$
and the graph of the initial velocity $y=u_0(s)$ have a common
linear segment $\bar \Omega = [s_1, s_2].$

Then at the moment $t=t_*$ at the point $x=x_*$ the component of the
density develops a $\delta$ -- singularity of amplitude
\begin{equation}\label{amplitude}
A(x_*)\, = \int\limits_{s_1}^{s_2}\,f(s)\,ds.
\end{equation}
\end{Theorem}

\proof

We are going to prove that
$$\int\limits_{\mathbb R}\rho_\sigma (t_*,x) \phi(x)\, dx\,\to \,A(x_*)
\phi(x_*),\qquad \sigma\to 0.$$

From (\ref{plotn})  we have %for every $\phi \in
%C_0^\infty(\mathbb R)$
$$\int\limits_{\mathbb R}\rho (t_*,x) \phi(x)\, dx\,=\,\dfrac1{\sqrt{2\pi t_*}\sigma}\,\int\limits_{\mathbb
R}\,\int\limits_{\mathbb{R}}f_0(s)e^{-\frac{(u_0(s)t_*+s-x)^2}{2\sigma^2t_*}}\phi(x)\,ds\,dx\,=$$$$
\dfrac1{\sqrt{2\pi t_*}\sigma}\,\int\limits_{\mathbb
R}\,\int\limits_{\bar\Omega}f_0(s)e^{-\frac{(u_0(s)t_*+s-x)^2}{2\sigma^2t_*}}\phi(x)\,ds\,dx\,+$$
$$\dfrac1{\sqrt{2\pi t_*}\sigma}\,\int\limits_{\mathbb
R}\,\int\limits_{\mathbb{R}\setminus\bar\Omega}f_0(s)e^{-\frac{(u_0(s)t_*+s-x)^2}{2\sigma^2t_*}}\phi(x)\,ds\,dx\,=
\,I_1\,+\,I_2.$$

First we analyze $I_1.$

$$I_1= \phi (x_*)\,\int\limits_{\bar\Omega}\,f_0(s)\,\int\limits_{\mathbb
R}\,\dfrac1{\sqrt{2\pi
t_*}\sigma}\,e^{-\frac{(u_0(s)t_*+s-x)^2}{2\sigma^2t_*}}\,dx\,ds\,+$$
$$
\dfrac1{\sqrt{2\pi
t_*}\sigma}\,\int\limits_{U_{x_*}(\varepsilon)}\,\int\limits_{\bar\Omega}f_0(s)
e^{-\frac{(u_0(s)t_*+s-x)^2}{2\sigma^2t_*}}\,ds\,(\phi(x)-\phi(x_*))\,dx\,+
$$
$$
\dfrac1{\sqrt{2\pi
t_*}\sigma}\,\int\limits_{U_{x_*}(\varepsilon)\mathbb R\setminus
U_{x_*}(\varepsilon)}\,\int\limits_{\bar\Omega}f_0(s)e^{-\frac{(u_0(s)t_*+s-x)^2}{2\sigma^2t_*}}\,ds\,(\phi(x)-\phi(x_*))\,dx\,=
\,I_{11}+I_{12}+I_{13},$$ where $U_{x_*}(\varepsilon)$ is an
$\varepsilon$ - neighborhood of $x_*.$ The integral $I_{11}$ is
equal to $A(x_*) \phi(x_*),$ with $A(x_*)$ specified in the
statement of Theorem \ref{T2}, $I_{12}$ and $I_{13}$ tend to zero as
$\sigma\to 0,$ as can be shown in a standard way.

Further,
$$I_2\,=\,\dfrac1{\sqrt{2\pi t_*}\sigma}\,\int\limits_{\mathbb
R}\,\int\limits_{-\infty}^{s_1}\,f_0(s)e^{-\frac{(u_0(s)t_*+s-x)^2}{2\sigma^2t_*}}\phi(x)\,ds\,dx\,+
$$
$$\dfrac1{\sqrt{2\pi
t_*}\sigma}\,\int\limits_{\mathbb
R}\,\int\limits_{s_2}^{\infty}\,f_0(s)e^{-\frac{(u_0(s)t_*+s-x)^2}{2\sigma^2t_*}}\phi(x)\,ds\,dx\,=
\,I_{21}\,+\,I_{22}.$$ Let us prove that $I_{2}$ vanishes as
$\sigma\to 0.$ Since $I_{21}$ and $I_{22}$ can be  analyzed
similarly, we consider only $I_{21}.$

For  $\,s\in U_{s_1-0}(\varepsilon), \,\varepsilon>0$ from
(\ref{U_eps}), $m=2,$  we have
$$
u_0(s)\,t_*+s-x\,=\,
\frac{1}{2}\,u^{(2)}_0(s_{**)}\,(s-s_*)^2\,t_*-(x-x_*),
$$
with $s_{**}\in U_{s_1-0}(\varepsilon).$

Thus,
$$
I_{21}=\, \dfrac1{\sqrt{2\pi
t_*}\sigma}\,f_0(s_1)\,\phi(x_*)\,\int\limits_{\mathbb
R}\,\int\limits_{-\infty}^{s_1}\,e^{-\frac{(\frac{u_0''(s_1-0)t_*}{2}\,(s-s_1)^2-(x-x_*))^2}{2\sigma^2t_*}}\,ds\,dx\,+
$$
$$
+\dfrac1{\sqrt{2\pi t_*}\sigma}\,\int\limits_{\mathbb
R}\,\int\limits_{-\infty}^{s_1}\,e^{-\frac{(\frac{u_0''(s_1-0)t_*}{2}\,(s-s_1)-(x-x_*))^2}{2\sigma^2t_*}}\,
(f_0(s)\,\phi(x)-f_0(s_1)\,\phi(x_*))\,ds\,dx\,+
$$
$$
+\dfrac1{\sqrt{2\pi t_*}\sigma}\,\int\limits_{\mathbb
R}\,\int\limits_{-\infty}^{s_1}\,\left(e^{-\frac{(u_0(s)t_*+s-x)^2}{2\sigma^2t_*}}\,-\,
e^{-\frac{(\frac{u_0''(s_1-0)t_*}{2}\,(s-s_1)-(x-x_*))^2}{2\sigma^2t_*}}\right)\,
f_0(s)\,\phi(x)\,ds\,dx=
$$
$$=\,I_{211}\,+\,I_{212}\,+\,I_{213}.
$$
To evaluate $I_{211}$ we use the formula $$\int\limits_{{\mathbb
R}}\,e^{-\frac{(a\,s^2\,-x)^2}{b^2}}\,ds\,=\,\frac{\sqrt{2}}{2}\,e^{-\frac{x^2}{2b^2}}\,\sqrt{\left|\frac{x}{a}\right|}
{\mathcal B}_{\frac{1}{4}}\left({-\frac{x^2}{2b^2}}\right),$$ where
$a\ne 0, b\ne 0$ are constants, ${\mathcal B}_{\frac{1}{4}} $ is the
Bessel function of the second type (\cite{Ryzhik}). Thus,
\begin{equation}\label{Mt*}
I_{211}\,=\,\sigma^{\frac{1}{2}}\,f_0(s_1)\,\phi(x_*)\,M(t_*,x_*),
\end{equation}
where
$$
M(t_*,x_*)\,=\,
 \frac{1}{2^{\frac{3}{4}}\pi^{\frac{1}{2}}}\frac{|u'(s_1)|^{\frac{1}{4}}}{|u''(s_1)|^{\frac{1}{2}}}\,L,\quad
 L\,=\,\int\limits_{\mathbb R}\,e^{-\frac{x^2}{2}}\,\sqrt{|x|}\,
{\mathcal B}_{\frac{1}{4}}\left({-\frac{x^2}{2}}\right)\,dx.$$ The
integral in the expression of $L$ converges, since the integrand is
finite at $x=0$ and decays exponentially at infinity.

The fact that $\,I_{212}\,$ and $\,I_{213}\,$ vanish as $\sigma\to
0$ can be proved routinely. $\square$

\medskip

\begin{remark} Under the assumptions of Theorem \ref{T1} one can similarly show (analogously to (\ref{Mt*})) that
$$
\rho_\sigma(t_*,x_*)\,\sim\,\sigma^{\frac{1}{m}}\,M(x_*,t_*,m)\,f_0(s_*)\,R_\sigma(t_*,x_*),\qquad\sigma\to
0,
$$
where $M(x_*,t_*,m)$ is a constant depending only on the properties
of $u_0$ and  $R_\sigma (t_*,x_*)$ tends to $\delta (x-x_*)$ as
$\sigma\to 0$ in ${\mathcal D}'({\mathbb R}).$ Thus, the $\delta$ -
function, arising from smooth initial data without linear segments,
has initially a zero amplitude.
\end{remark}

%$$=\,\dfrac1{(\sqrt{2\pi
%t_*}\sigma)}\,\int\limits_{U_{x_*}(\varepsilon)}
%\,\int\limits_{U_{s_*}(\varepsilon)}\,f_0(s)\,e^{-\frac{(u''_0(s_{**}))^2\,t_*}{8\,\sigma^2}\,(s-s_*)^4}\,\phi(x)\,ds\,dx+
%o(\sigma).
%$$

%\subsection{Riemann problem with nonconstant states}

\section{The Hugoniot conditions and the spurious
pressure}\label{spurpress} As follows from the results of Sec.2, if
$f_{FP}$ and $u_{FP}$ are smooth, they solve the pressureless gas
dynamics system. Now we ask the question which system satisfies the
FP-generalized solution  with jumps obtained in Sec.4.

 The
system of conservation laws (\ref{sist_pred1}), (\ref{sist_pred2})
implies two Hugoniot conditions that should be held on the jumps of
the solution \cite{Rogd}. This signifies that the solution satisfies
the system in the sense of integral identities. If we denote by
$\mathfrak D$ the velocity of the jump and by $[h(y)]=h(y+0)-h(y-0)$
the value of the jump, then the continuity equation and the momentum
conservation give $[f]\,{\mathfrak D}\,=\,[fu]$ and
$[fu]\,{\mathfrak D}\,=[fu^2],$ respectively.

In the case $u_2>0$ the velocity is continuous, therefore the
Hugoniot conditions hold trivially.

We should check these conditions for the jumps in the case $u_2<0.$
An easy computation shows that the first one is satisfied:  for the
point $\hat{x}_1=u_1t$ we have:

$${\mathfrak D}\,=\,\dfrac{(f_1+f_2)(u_1+u_2)-(2f_1+f_2)u_1-(f_1+f_2)u_2}{-f_1}=u_1,$$

and for $\hat{x}_2=(u_1+u_2)t$:

$${\mathfrak D}\,=\,\dfrac{(2f_1+f_2)u_1+(f_1+f_2)u_2-f_1u_1}{f_1+f_2}=u_1+u_2.$$

However, the second Hugoniot condition does not hold. To understand
the reason for this let us estimate the integral term in
(\ref{sist_obw2})  in the case $u_2<0$ as $\sigma \to 0$:

$$
\int\limits_{\mathbb R}\,(u-\hat u(t,x))^2\,P_x(t,x,u)\,du=$$
$$
\frac{1}{\sqrt{2\pi t}\sigma}\,\int\limits_{\mathbb
R}\,f_0(s)(u_0(s)-\hat
u(t,x))^2\,\left(e^{-\frac{(u_0(s)t+s-x)^2}{2\sigma^2t}}\right)_x\,ds=
$$
$$
 -\,\frac{1}{\sqrt{2\pi t}\sigma}\,\int\limits_{\mathbb
R}\,f_0(s)((u_0(s)-u_{FP}(t,x)) +(u_{FP}(t,x)-\hat
u(t,x)))^2\,\left(e^{-\frac{(u_0(s)t+s-x)^2}{2\sigma^2t}}\right)_s\,ds=
$$
$$
 -\,\frac{1}{\sqrt{2\pi t}\sigma}\,\int\limits_{\mathbb
R}\,f_0(s)(u_0(s)-u_{FP}(t,x))^2\,\left(e^{-\frac{(u_0(s)t+s-x)^2}{2\sigma^2t}}\right)_s\,ds-
$$
$$
2\,(u_{FP}(t,x)-\hat u(t,x))\,\frac{1}{\sqrt{2\pi
t}\sigma}\,\int\limits_{\mathbb
R}\,f_0(s)(u_0(s)-u_{FP}(t,x))\,\left(e^{-\frac{(u_0(s)t+s-x)^2}{2\sigma^2t}}\right)_s\,ds+
$$
$$
\,(u_{FP}(t,x)-\hat u(t,x))^2\,\frac{1}{\sqrt{2\pi
t}\sigma}\,\int\limits_{\mathbb
R}\,f_0(s)\,\left(e^{-\frac{(u_0(s)t+s-x)^2}{2\sigma^2t}}\right)_s\,ds=
$$
$$I_1+I_2+I_3.$$
The integrals $I_2$ and $I_3$ tend to zero as $\sigma\to 0$ due to
properties of the Riemann data since $\hat u(t,x)\to u_{FP}(t,x)$
for almost all $x\in \mathbb R.$ Let us estimate $I_1.$
$$I_1\,
=-\,\frac{1}{\sqrt{2\pi t}\sigma}\,\int\limits_{u_2 t}^0\,f_1\,(u_1-
u_{FP}(t,x))^2\,\left(e^{-\frac{(u_1\,t+s-x)^2}{2\sigma^2t}}\right)_s\,ds\,-$$$$
\,\frac{1}{\sqrt{2\pi t}\sigma}\,\int\limits_{0}^{-u_2
t}\,(f_1+f_2)\,((u_1+u_2)-
u_{FP}(t,x))^2\,\left(e^{-\frac{((u_1+u_2)\,t+s-x)^2}{2\sigma^2t}}\right)_s\,ds
=
$$
$$
-\,\frac{1}{\sqrt{2\pi t}\sigma}\,\frac{f_1(f_1+f_2)^2
u_2^2}{(2f_1+f_2)^2}\left(e^{-\frac{(u_1\,t-x)^2}{2\sigma^2t}}-e^{-\frac{((u_1+u_2)\,t-x)^2}{2\sigma^2t}}\right)-$$
$$\,\frac{1}{\sqrt{2\pi t}\sigma}\,
\frac{f_1^2(f_1+f_2)
u_2^2}{(2f_1+f_2)^2}\left(e^{-\frac{(u_1\,t-x)^2}{2\sigma^2t}}-e^{-\frac{((u_1+u_2)\,t-x)^2}{2\sigma^2t}}\right).$$
Thus,
$$I_1\,\to\,-\frac{f_1(f_1+f_2)
u_2^2}{(2f_1+f_2)}(\delta(x-(u_1+u_2)\,t)\,- \, \delta(x-u_1\,t)
),\quad \sigma\to 0,$$ in the distributional sense.

Thus, the integral term corresponds to a spurious pressure $p(t,x)$
between the  jumps $x=(u_1+u_2)\,t$  and $x=u_1\,t,$ namely,
\begin{equation}\label{press}p(t,x)\,=\,\frac{f_1(f_1+f_2)
u_2^2}{(2f_1+f_2)}(\theta(x-(u_1+u_2)\,t)\,- \, \theta(x-u_1\,t)),
\end{equation}
see Fig.\ref{pressure1}.

\begin{figure}[h]
\centerline{\includegraphics[width=1 \columnwidth]{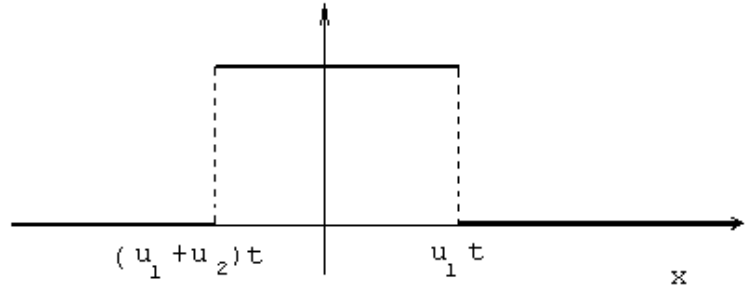}}
\caption{Spurious pressure, $u_2<0.$}\label{pressure1}
\end{figure}
 The Hugoniot condition $[fu]\,{\mathfrak D}\,=\,[\,fu^2+p\,]$ is
satisfied with this kind of pressure.

Thus, we get the following theorem.

\begin{Theorem} The generalized FP-solution to the Riemann problem  (\ref{f0_razrG}), (\ref{u0_razrG}) with constant left-hand
and right-hand initial states  for the pressureless gas dynamics
system in the case of a discontinuous velocity ($u_2<0$) solves in
fact the gas dynamics system with a pressure defined by
(\ref{press}).
\end{Theorem}

The analogous calculations for the case of rarefaction $u_2<0$ show
that $$I_1=-\,\frac{1}{\sqrt{2\pi
t}\sigma}\,\int\limits_{-\infty}^0\,f_1\,(u_1-
u_{FP}(t,x))^2\,\left(e^{-\frac{(u_1\,t+s-x)^2}{2\sigma^2t}}\right)_s\,ds\,-$$$$
\,\frac{1}{\sqrt{2\pi
t}\sigma}\,\int\limits_{0}^{+\infty}\,(f_1+f_2)\,((u_1+u_2)-
u_{FP}(t,x))^2\,\left(e^{-\frac{((u_1+u_2)\,t+s-x)^2}{2\sigma^2t}}\right)_s\,ds=
$$
$$\,\frac{1}{\sqrt{2\pi
t}\sigma}\,\left((u_1-\frac{x}{t})^2
\,f_1\,e^{-\frac{(u_1\,t-x)^2}{2\sigma^2t}}-
((u_1+u_2)-\frac{x}{t})^2
\,f_1\,e^{-\frac{((u_1+u_2)\,t-x)^2}{2\sigma^2t}}\right).$$ Here we
use the FP-solution $(f_{FP}, u_{FP}),$ obtained in Sec.4.

Thus,  $I_1\to 0$ as $\sigma\to 0,$ therefore the integral term
vanishes in the case $u_2<0.$

\section{Sticky particles model vs noninteracting particles}

In our model the particles are allowed to go through the
discontinuity as one particle does not feel the others. However, in
the frame of the sticky particles model  the particles meeting each
other are assumed to stick together on the jump \cite{ERykovSinai}.
The noninteracting particles model and the sticky particles model
are equivalent to the same system (\ref{sist_pred1}),
(\ref{sist_pred2}) for smooth densities and velocities, however, if
the initial data have a jump, the solutions behavior differs
drastically.
\subsection{Riemann problem with constant states}
Nevertheless, we can study the solution to the Riemann problem in
the case of $u_2<0$ for the sticky particles model, too, using  the
solution obtained above for the non-interacting model. Indeed, the
jump position $x_j(t)$ is  a point between $x_1(t)=(u_1+u_2)\,t$ and
$x_2(t)=u_1\,t.$ The mass $m(t)$ accumulates in the jump due to the
impenetrability of the discontinuity as
$$m(t)=(x_j(t)-(u_1+u_2)t)(f_1+f_2)+$$

$$+(u_1t-x_j(t))f_1\,+\, f_3\,=-((u_1+u_2)(f_1+f_2)-u_1\,f_1)\,t+$$

\begin{equation}\label{m(t)}+x(t) f_2\,+\, f_3\,=\, -[uf]\,t\,+\,[f]\,x_j(t)\,+\,f_3,
\end{equation} where
$[\,\,]$ stands for the jump value, and $m(0)=f_3,\,x_j(0)=0.$
%Therefore
%due to the initial conditions $x(0)=0,\,$ we have
%$$
%m(t)=-[uf]t+[f]x(t)+m(0)
%$$
Further, if we change heuristically the overlapped  mass between
$x_1$ and $x_2$ to the mass concentrated at a point (see
Fig.\ref{sticky}), then from the condition of equality of momenta in
both cases we can find the velocity of the point singularity:
$$(u_1+u_2)(f_1+f_2)(x_j(t)-(u_1+u_2)t)+$$
$$+u_1\,f_1\,(u_1t-x_j(t))\,=$$
$$=-[u^2 f]\,t\,+\,[u f]\, x_j(t)\,=\,m(t)\,\dot x_j(t).$$
%$$=D\left(-((u_1+u_2)(f_1+f_2)-u_1\,f_1)+Df_2\right).$$
Thus, to find the position of the point singularity we get the
equation
$$
([f]\,x_j(t)\,-\,[uf]\,t \,+\, f_3\,)\,\dot x_j(t)\,=\,[u f]\,
x_j(t)\,-\,[u^2 f]\,t,
$$
subject to the initial data $x_j(0)=0.$ The respective solution is
\begin{equation}\label{x(t)[f]}x_j(t)=\frac{1}{[f]}\,\left([uf]\,t\,-f_3\,
+\,\sqrt{f_3^2\,-\,2\,[uf]\,f_3\,t\,+\,([uf]^2-[f][u^2f])\,t^2}\right),\,
\mbox{if}\,\,[f]\ne 0,
\end{equation}
and
\begin{equation}\label{x(t)[0]}
x_j(t)=\frac{[u^2]\,f \,t^2}{2([u]\,f\,t\,-\,f_3)},\quad
\mbox{if}\,\,[f]\,=\,0.\end{equation} In particular,  from the
latter formula in the case $f_3=0$ we get the known expression for
the velocity of the jump \cite{Rogd}:
$$\dot x_j(t)\,=\,\frac{2u_1+u_2}{2}\,=\,\frac{u_{left}+u_{right}}{2}.$$
%$$x(t)=\frac{1}{[f]}\,\left([uf]\pm\sqrt{[uf]^2-[f][u^2f]}\right)\,t,$$
It can be checked that  $x_1(t)<x_j(t)<x_2(t).$ The condition
expressed by these inequalities is equivalent to the Lax stability
condition $u_1<\dot x_j(t)<u_1+u_2.$

\begin{figure}[h]
\centerline{\includegraphics[width=1\columnwidth]{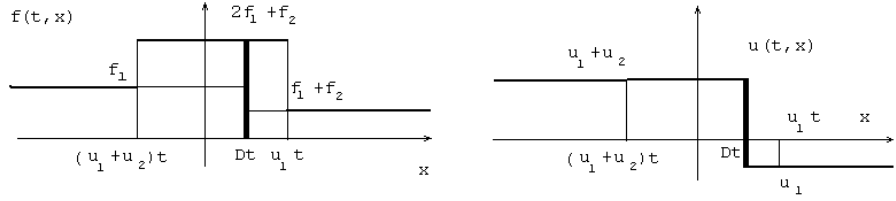}}
\caption{Changing the noninteracting particles model to the sticky
particles model.} \label{sticky}
\end{figure}

The formulas describing  the amplitude of the delta-function in the
density component and the singularity position  obtained earlier in
\cite{huang},\cite{shelkv},\cite{danilov} give the same result.
Moreover, our method allows to find the jump position in a unique
way (in contrast to the method used in \cite{danilov}).
% since the value of $\dot x (0)$
%is not a fr

It is worth mentioning that the spurious pressure (\ref{press}) does
not arise in the sticky particles model.

\begin{remark} The case of the singular Riemann problem  with $f_3 \ne 0$ is in fact more different
from the case of the regular Riemann problem than it seems at first
glance.  Let us begin with the constant initial density ($[f]=0$),
when the trajectory satisfies (\ref{x(t)[0]}). It is natural to set
the amplitude of the initial $\delta$ -- function $f_3$ greater or
equal than zero. Then the trajectory $x(t)$ is continuous since the
denominator in (\ref{x(t)[0]}) does not vanish. However, if we take
$f_3<0,$ then the trajectory goes to infinity at the finite moment
where $m(t)$ vanishes, and then the trajectory jumps to infinity of
the other sign.

Further, if $[f]\ne 0,$ then we have to use the formula
(\ref{x(t)[f]}). It can be checked that it is possible to find
values $u_1, u_2, f_1, f_2$ (for example, $u_1=-1, u_2=-2, f_1=2,
f_2=-1.8 $) such that the expression under the square root vanishes
within a finite time $t_*.$ However, as follows from (\ref{m(t)}),
(\ref{x(t)[f]}),
\begin{equation}\label{m(t)[f]}m(t)\,=\,\sqrt{f_3^2\,-\,2\,[uf]\,f_3\,t\,+\,([uf]^2-[f][u^2f])\,t^2},\,
\end{equation}
therefore within the same time $t_*$ the amplitude of the $\delta$
-- function becomes zero. Thus, at the moment $t_*$ we have to set a
new Riemann problem with the jump at the point $x(t_*).$

%Let us note that in  formula (\ref{x(t)[f]}) the coefficient
%$$([uf]^2-[f][u^2f])=
%f_{left}\,f_{right}\,(u_{right}\,-\,u_{left})^2$$ is positive  in
%the absence of vacuum in initial state.
\end{remark}
\subsection{Riemann problem with non-constant states}
Now we extend formulas (\ref{x(t)[f]}) and (\ref{x(t)[0]}) to the
case of the Riemann problem with non-constant left and right states:
\begin{equation}
f_0(x)=f_1(x)+\theta(x) f_2(x)+ f_3\delta(x),
\end{equation}
\begin{equation}
u_0(x)=u_1(x)+\theta(x) u_2(x),
\end{equation}
where $u_1(x),$ $u_2(x),$ $f_1(x),$ $f_2(x)$ are smooth
%continuously
%differentiable
functions, $f_3$ is a real constant. We restrict ourselves to a
situation that is quite similar to the case of constant states.
Namely, we assume that for every $x\in\mathbb R$ and $t>0$ the
straight line $y=\frac{x-s}{t}$ has at most two common points with
the graph of the function $y=u_0(s),$ moreover, let
$u_2^0=\lim\limits_{x\to 0+} u_2(x)<0$ and assume that the
intersection points $s_-(t,x)$ and $s_+(t,x)$ lie on either side of
the origin $x=0$ (see Fig.\ref{RPgen}).

\begin{figure}[h]
\centerline{\includegraphics[width=1\columnwidth]{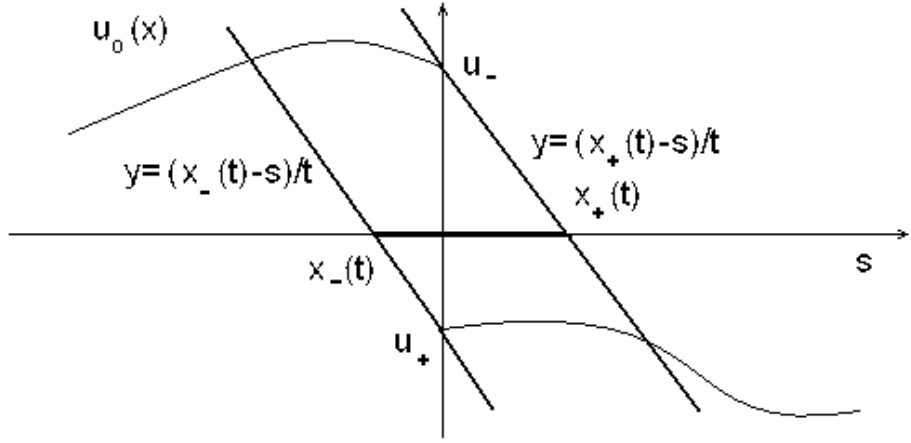}}
\caption{The overlapping domain in the Riemann problem}
\label{RPgen}
\end{figure}

So, every point $x$ at the moment $t>0$ lies in the overlapping
domain
$$\,D\,=\,[x_-(t), x_+(t)],\quad x_-(t)=u_+ t, \, x_+(t)=u_- t,$$ if
the line $y=\frac{x-s}{t}$ and the graph of $y=u_0(s)$ have two
points of intersection, $s_-(t,x)<0$ and $s_+(t,x)>0.$ We set
$u_-\,\equiv \,\lim\limits_{x\to 0-}\,u_1(x)$ and
$u_+\,\equiv\,u_-+u_2^0.$ Further, we assume that for any fixed $t$
the number $K$ of points $\bar x_k(t)\in D $ such that the straight
line and the graph of the initial velocity have a common linear
segment is finite.
%At last, we assume that the domain $D$ does not meet  singularities
 Let us again denote by $x_j(t)$ the position of the singularity that
 should change the overlapping domain $D$ in the sticky particles
 model.
 Then the conservation of mass gives
\begin{equation}\label{m(t)RPnonconst}m(t)\,=\,\int\limits_{u_+ t}^{x_j(t)}\,f(s_+(t,x))\,dx
\,+\,\int\limits_{x_j(t)}^{u_- t}\,f(s_-(t,x))\,dx\,+\,
\sum\limits_{k=1}^{K}\,A_k(t,\bar x_k(t))\,+\,f_3,
\end{equation}
where $\,A_k(t,\bar x_k(t))\,$ is the amplitude of the $\delta$ --
function formed at the point $\bar x_k(t),$ calculated by using the
formula (\ref{amplitude}).

Further, from the conservation of momentum we have
\begin{eqnarray}\label{momentRPnonconst}m(t)\,\dot x_j(t)\,=\,\int\limits_{u_+ t}^{x_j(t)}\,f(s_+(t,x))\,u(s_+(t,x))\,dx
\,+\,\int\limits_{x_j(t)}^{u_-
t}\,f(s_-(t,x))\,u(s_-(t,x))\,dx\,+\,\end{eqnarray}
$$+\,\sum\limits_{k=1}^{K}\,A_k(t,\bar x_k(t))\,\dfrac{d}{dt}\bar
x_k(t),$$ moreover, the velocity $\dfrac{d}{dt}\bar x_k(t)$ can be
found by the formula
$$\dfrac{d}{dt}\bar x_k(t)\,=\,u_0(s_*(t,\bar x_k(t))),$$
where $s_*(t,\bar x_k(t))$ is the coordinate of the point where the
graph of $y=u_0(s)$ and the line $y=\frac{x-s}{t}$ have a common
linear segment. If the graph of $y=u_0(s)$ does not contain any
linear segments, then the respective parts in formulas
(\ref{m(t)RPnonconst}) and (\ref{momentRPnonconst}) vanish. Thus, it
is sufficient to substitute (\ref{m(t)RPnonconst}) into
(\ref{momentRPnonconst}) to get  the integro-differential equation
that governs the singularity position. This equation should be
considered together with the initial condition $x_j(0)=0.$

\subsection{Evolution of the singularity formed from smooth data}

As we have seen in Sec.{\ref{Arising}}, if at a point $x$ and a
moment of time $t_*$  starting from smooth initial data loses its
smoothness, then there arises a gradient catastrophe in the velocity
component (the derivative becomes unbounded), whereas in the density
component there arises a $\delta$ -- singularity.  In the framework
of the pressureless gas dynamics for $t>t_*$ the $\delta$ --
singularities encompass the overlapping domain $D$ and  in the
overlapping domain the spurious pressure (given by the integral term
discussed in Sec.\ref{spurpress}) appears. In the sticky particles
model we have to collapse the overlapping domain to the one point,
where the whole mass of $D$ accumulates. The position of this new
singularity should be found based on the conservation of mass and
momentum.

For the sake of simplicity we assume that every straight line
$y=\frac{x-s}{t}$ intersects the graph of the smooth initial
velocity $y=u_0(s)$ at most three times. For every fixed $x$
initially the intersection point is unique ($0<t<t_*$), then at the
moment $t=t_*$ the straight line becomes  tangent to the graph of
the initial velocity in a certain point, and for $t>t_*$ we have
three intersection points. Our aim is to find the position of the
singularity and the amplitude of the $\delta$ -- function in the
density component.

We denote as before by $x_-(t)$ and $x_+(t)$ the endpoints of the
domain $D, $ $x_j(t)$ -- the position of the new singularity,
$A(x_-(t))$ and $A(x_+(t))$ -- the amplitude of $\delta$-- functions
produced at the point where the graphs of $y=u_0(s)$ and
$y=\frac{x_\pm(t) -s}{t}$ have a common linear segment (for fixed
$t\ge t_*$). Further, let $s_-(t,x),$ $s_0(t,x),$ $s_+(t,x)$ be the
subsequent point of intersection ($ s_- < s_0 < s_+.$) Let $m(t)$ be
the amplitude of the $\delta$ - function in the density component.

Then, the conservation of mass gives
\begin{equation}\label{sm_dat_m(t)}
m(t)\,=\,A(x_-(t))\,+\,A(x_+(t))\,+
\end{equation}
$$\,\int_{x_-(t)}^{x_j(t)}\,(f(s_+(t,x))\,+\,f(s_0(t,x))\,dx \,+\,
\int_{x_j(t)}^{x_+(t)}\,(f(s_-(t,x))\,+\,f(s_0(t,x))\,dx,$$
 where $A(x_\pm (t))$ ca be found by formula
(\ref{amplitude}).

From the conservation of momentum analogously to
(\ref{momentRPnonconst}) we have
\begin{eqnarray}\label{moment_smooth_data}m(t)\,\dot x_j(t)\,=\,\int\limits_{x_-(t)}^{x_j(t)}\,
\big(f(s_+(t,x))\,u(s_+(t,x))+f(s_0(t,x))\,u(s_0(t,x))\big)\,dx
\,+\end{eqnarray}$$\,\int\limits_{x_j(t)}^{x_+(t)}\,\big(f(s_-(t,x))\,u(s_-(t,x))+f(s_0(t,x))\,u(s_0(t,x))\big)\,dx\,+\,$$
$$\,A(x_-(t))\, u_0(s_-(t,x_-(t)))\,       +\,A(x_+(t))\, u_0(s_+(t,x_+(t))).$$

Thus,  equations (\ref{sm_dat_m(t)}),  (\ref{moment_smooth_data})
and the initial conditions $m(t_*)=A(x_*(t_*)),$ $x_j(t_*)=x_*$
 define the position  and
the amplitude of the singularity of the $\delta$ - function in the
component of density. Here  $x_*$ is  a point such that
$y=\frac{x_*-s}{t_*}$ and $y=u_0(s)$ have a common point $s_*$ or a
common linear segment $[s_-,s_+]$ such that the derivatives of both
functions are equal at $s_*$ or on $[s_-,s_+]$ , $A(x_*)$ is defined
in the statement of Theorem \ref{T2}.

\begin{remark}
If we want to consider the global evolution of the solution to the
Burgers equation itself, we should set $f_0=\rm const.$ The
continuity equation plays here an auxiliary role. We are not
interested in the properties of the density $f,$ which is constant
everywhere except for domains of vacuum and except for the points
giving  rise to the $\delta$ - singularity.
\end{remark}

\section{Extension to more general scalar conservation law}
Let us consider the following equation
\begin{equation}\label{ext_Burgers}
\partial_t v+(G(v),\nabla)v=0,
\end{equation}
subject to initial data  $v(x,0)=v_0(x),$  where
$v(x,t)=(v_1,...,v_n)$ is a vector-function
$\mathbb{R}^{n+1}\rightarrow\mathbb{R}^n,$ $G(v)$ is a
non-degenerate differential  mapping from ${\mathbb R}^n\,$ to
$\,{\mathbb R}^n, $ such that its Jacobian satisfies ${\rm det}
\frac{\partial G_i(v)}{\partial v_j} \ne 0,$ $i,j=1,...,n.$
%where $(X(t),U(t))$ runs in the phase space
%$\mathbb{R}^n\times\mathbb{R}^n,$  $\sigma$ is a real strictly
%positive constant and $(W)_t=(W_k)_t$, $k=1,...,n$ is an n -
%dimensional Brownian motion, $\,t>0.$

We can multiply (\ref{ext_Burgers}) by $\nabla G_i(v),\,i=1,...,n,$
to get
\begin{equation}\label{ext_Burgers1}
\partial_t G(v)+(G(v),\nabla)G(v)=0.
\end{equation}
Thus, we can introduce a new vectorial variable $u=G(v)$ to reduce
the Cauchy problem for (\ref{ext_Burgers1}) to (\ref{equ_Burg}) with
$u_0(x)=G(v_0(x)).$ The stochastic perturbation for
(\ref{ext_Burgers1}) is (\ref{SDU}) with  $U$ replaced by $G(V).$
%it
%is equivalent to (\ref{SDUExt}).

Therefore we  find the representation of the solution to the
stochastically perturbed along the characteristics of equation
(\ref{ext_Burgers1}) using the formula (\ref{sol_u_sdu}) with
$G(v_0(x))$ instead of $u_0(x).$

Thus, we can apply the results obtained in the previous sections to
the investigation of the Riemann problem and the arising of
singularities for the following analogue of the pressureless gas
dynamics system:
\begin{equation}\label{newpressureless}
g_t+ {\rm div}_x (g G(v))=0,\quad (g G(v))_t\,+\,\nabla_x (g
\,v\,\otimes\,G(v))=0,
\end{equation}
where $g(t,x)$ is a scalar function that can be interpreted as a
density.

Thus, just as we relate with the non-viscous Burgers equation the
system of pressureless gas dynamics, so also with the equation
(\ref{ext_Burgers}) one can relate the system
(\ref{newpressureless}).

To obtain the  solution to the Cauchy problem for
(\ref{ext_Burgers}) itself we have to consider the solution to the
Cauchy problem for (\ref{newpressureless})  with the data
$(g_0,\,G(v_0) ),\,g_0 \equiv {\rm const}$ and then  perform the
inverse transform $v(t,x)=G^{-1}(u(t,x)).$
%Finally let us mention that setting $\,g_0(x)=\rm const\,$ we can
%investigate both the Riemann problem and the singularity formation
%for (\ref{ext_Burgers}).

\section{Conclusion}

1. Let us notice that the solution of the Riemann problem
%withconstant states
 for the
pressureless gas dynamics system obtained in this paper for the  1D
case satisfies the entropy condition
\begin{equation}\label{E-cond}
\frac{u(t,x_2)-u(t,x_1)}{x_2-x_1}\,\le\,\frac{1}{t},
\end{equation}
for any sufficiently small $x_1$ and $x_2$ (e.g.\cite{huang}) and
the balance relations on the jump that the definition of solution in
the sense of integral identity implies (\cite{shelkv}) are satisfied
as well. It is known that these conditions do not guarantee
uniqueness (\cite{huang}, \cite{danilov}).  However, our solution is
unique both in the case of rarefaction and compression. In the case
of rarefaction it is automatically self-similar (we recall that the
assumption of the self-similarity implies uniqueness
\cite{ShengZhang}). In the case of contraction the problem of
uniqueness was open for the solution to the singular Riemann
problem, where  a non-zero mass is concentrated on the jump at the
initial time. As was noticed in \cite{danilov}, for the uniqueness
one must prescribe the derivative of the amplitude of the $\delta$ -
function. In our framework the solution is unique and the value of
the derivative of the amplitude of the $\delta$ - function follows
from the expression for the amplitude itself.

\medskip

2. In \cite{Derm} an analog of the system (\ref{sist_obw1}),
(\ref{sist_obw2}) (without the integral term) in the 1D case was
obtained. Namely, it was proved that for { smooth} initial data a
local in time strong solution $(\rho(t,x), \,u(t,x))$ to this system
can be constructed by means of a nonlinear diffusion process
$$X_t=X_0+\int\limits_0^t\,\mathbb E[(u_0(X_0))|X_s]\,ds\,+\,\sigma
W_t,\quad \mathcal L(X_0)=\rho(0,x),$$ so that $\rho(t,x)$  is the
probability density of the diffusion process $X_t$ and $
u(t,x)=E[(u_0(X_0))|X_t=x]$ (with $E[ .,| \,..]$ standing for
conditional expectation and $\,{\mathcal L}(X_0)\,$ for the
probability density of $X_0$). In fact, this result relates to our
Proposition 2, since we have shown that the integral term arises
only for discontinuous data.

Further, in \cite{DermPD} the system (\ref{sist_obw1}),
(\ref{sist_obw2}) (without the integral term) was considered in any
dimension. In this work there was constructed a global weak solution
using  discrete approximations, and the interaction of particles is
given by a sticky particles dynamics.

3. There exist formalisms to represent solutions of parabolic PDE's
as the expected value of functionals of stochastic processes (see
e.g. \cite{Freidlin}, \cite{Freidlin1}, \cite{Belopolskaya},
\cite{Belopolskaya99} and references therein). In particular,
 in \cite{Iyer} one can find a recent result  concerning the stochastic formulation of the
viscous Burgers equation. An alternative approach to the stochastic
formulation for a much more wider class of parabolic equations and
systems can be found in  \cite{AlbeverioBelopolskaya}.

4. We would also like to mention the paper \cite{CKR}, where a
numerical method of particles for the solution of the pressureless
gas dynamics in 1D an 2D case has been developed. The method is
mostly inspired by \cite{ERykovSinai} and in fact in this paper the
problem of transition from the non-interacting particles to the
sticky particle model was solved numerically.

{\bf Acknowledgements}  We are grateful to Ya.Belopolskaya,
V.Danilov, A.Kurganov, E.Panov and V.Shelkovich for fruitful
discussions.

\end{document}